\documentclass[12pt]{amsart}
\usepackage{amssymb}
\usepackage[latin1]{inputenc}
\usepackage{amsthm}
\usepackage{amsfonts}
\usepackage{mathrsfs}
\usepackage{graphicx}
\usepackage{amsmath}
\let\mathcal\mathscr
\def\bB{{\mathbb B}}
\def\bC{{\mathbb C}}
\def\bG{{\mathbb G}}

\def\bP{{\mathbb P}}
\def\bN{{\mathbb N}}
\newtheorem{thm}{Theorem}[section]

\def\Z{{\bf Z}}
\def\N{{\bf N}}

\def\C{{\mathbb C}}

\def\codim{\mathop{\rm codim}\nolimits}

\def\Im{\mathop{\rm Im}\nolimits}

\def\Jac{\mathop{\rm Jac}\nolimits}

\def\lra{\longrightarrow}

\def\Pic{\mathop{\rm Pic}\nolimits}

\def\Sing{\mathop{\rm Sing}\nolimits}
\def\Sym{\mathop{\rm Sym}\nolimits}

\def\tilde{\widetilde}
\def\eps{\varepsilon}
\def\phi{\varphi}
\def\a{{\alpha}}

\def\i{{\infty}}

\def\0{{\bf 0}}
\def\t{{\bf t}}

\def\cA{{\mathcal A}}

\def\cD{{\mathcal D}}

\def\cO{{\mathcal O}}

\def\cX{{\mathcal X}}
\def\cY{{\mathcal Y}}

\def\Base{\mathop{\rm Base}\nolimits}

\def\rank{\mathop{\rm rank}\nolimits}

\def\ge{\geqslant}

\numberwithin{equation}{section}

\newtheorem{theorem}[thm]{Theorem}

\newtheorem{conjecture}[thm]{Conjecture}

\newtheorem{cor}[thm]{Corollary}

\newtheorem{definition}[thm]{Definition}
\newtheorem{example}[thm]{Example}

\newtheorem{lemma}[thm]{Lemma}

\newtheorem{proposition}[thm]{Proposition}
\newtheorem{prop}[thm]{Proposition}
\newtheorem{rmk}[thm]{Remark}
\newtheorem{remark}[thm]{Remark}

\newtheorem{pb}[thm]{Problem}

\begin{document}
\thanks{{\it Mathematics Subject Classification (2000)}: Primary: 14J70,
32Q45} \keywords{Hyperbolicity of projective varieties; holomorphic mappings; jet bundles; Green-Griffiths conjecture, Kobayashi conjecture.}
\title [Generalized Demailly-Semple jet bundles and holomorphic mappings]{Generalized Demailly-Semple jet bundles and holomorphic mappings into complex manifolds}
\author{Gianluca Pacienza and Erwan Rousseau}
\date{}
\maketitle

\begin{abstract}
Motivated by the Green-Griffiths conjecture,
we study maximal rank holomorphic maps from $\C^p$ into complex manifolds. When $p>1$ such maps should in principle be more tractable than entire curves. We extend to this setting the jet-bundles techniques introduced by Semple, Green-Griffiths and Demailly. Our main application is the non-existence of maximal rank holomorphic maps from $\C^2$ into the very general degree $d$ hypersurface in $\bP^4$, as soon as $d\geq 93.$   
\end{abstract}

\tableofcontents

\section{Introduction}\label{S:intro}
%
The motivation of this paper is the study of holomorphic maps $ f : \bC^{p} \rightarrow X $, $1\leq p\leq n$, where $X$ is a complex manifold of dimension $n$. Such a study is of interest for important problems related to hyperbolicity-type properties of complex manifolds. If $X$ is of general type, i.e. its canonical bundle $K_{X}$ is big, one can hope following Green and Griffiths that the following is true:
\begin{conjecture}
Let  $ f : \bC^{p} \rightarrow X $ be a holomorphic map in a compact complex manifold of general type. Then the image of $f$ lies in a proper subvariety $Y_{f}$ of $X$.
\end{conjecture}
The equidimensional  case $p=n$  has been settled by Kobayashi-Ochiai in \cite{KO75}. 
In the litterature, many articles are concerned with the case $p=1$ of entire curves, which is the most difficult one since by the following stronger conjecture it implies the other cases.
\begin{conjecture}
Let  $X$ be a compact complex manifold of general type. Then there exists a proper subvariety $Y \subset X$ such that for every holomorphic map $ f : \bC \rightarrow X $ the image of $f$ lies in $Y$.
\end{conjecture}
Despite the efforts of several people, this conjecture seems out of reach for the moment. We refer the reader to the bibliography of Demailly's  notes \cite{De95} for the precise references. 
Another important problem is the study of Kobayashi-Eisenman $p$-measure hyperbolicity. Generalizing the construction of the Kobayashi pseudometric, Eisenman introduced on each complex manifold of dimension $n$ a $p$-dimensional pseudometric for any $1\leq p \leq n$ (see \cite{E}). While a substantial work has been done on the case $p=1$ of the Kobayashi pseudometric, it seems that much less is known on the Eisenman pseudometrics.

The goal of our work is twofold. On the one hand we generalize some existing techniques to deal with the intermediate cases $2\leq p \leq n-1$, which in principle should be more tractable. Then we put these general constructions at work to extend to arbitrary $p$ several results which were obtained in the case of entire curves (e.g. \cite[Theorem 9.1]{De95}, \cite[Theorem 1]{Div}, \cite{DPP}, \cite[Theorem 6]{PR}). 
Our main application is the study of maps from $\bC^2$ to very general hypersurfaces in $\bP^4$.

More precisely, in the first part of this paper we extend to our more general setting the jet differentials techniques described in \cite{De95} to study entire curves. Namely we give the construction of the jet bundles $X_{k}^p$ as a tower of Grassmannians over $X$ and of the jet differentials bundles $E_{p,k,m}$ of invariant algebraic differential operators of order $k$ and degree $m$ acting on germs of holomorphic maps from $\bC^p$ to $X$ (see \S2 for the precise definitions). The construction of the jet bundles $X_{k}^p$ already appeared in a different setting in \cite{ASS}. The main technical difference with respect to the case $p=1$, extensively  studied in \cite{De95}, is the appearance of the effective locus $Z_{k}^p \subset X_{k}^p$, which is the closure of the set of points corresponding to the liftings of points of $X$ by holomorphic mappings.

Our main general results may be summarized in the following statement.

\begin{theorem}\label{thm:gen}
Let $X$ be a complex compact manifold, $A$ an ample line bundle on $X$, $f:\C^p\to X$ a maximal rank holomorphic mapping, i.e. whose differential is of maximal rank at some point of $X$, and $1 \leq p\leq \dim X, k$ and $m$ positive integers. Let $\pi_{0,k}:  X_{k}^p \rightarrow X$ be the natural projection.  Then for every non-zero global section
$$
 P\in  H^0(Z_{k}^p, \mathcal{O}_{Z_{k}}(m)\otimes \pi_{0,k}^*A^{-1})\subset H^0(X,E_{p,k,m}\otimes A^{-1}),
$$
the image of the $k$-th lifting $f_{[k]}$ of $f$ lies in the zero set of $P$. In other words, $P$ gives a non-trivial algebraic differential equation for maximal rank holomorphic mappings $f:\C^p\to X$.
\end{theorem}

We believe that these general constructions may lead to effective results for classes of varieties for which the study of entire curves happens to be very delicate.
For our main application, we use  
 these tools  here to obtain a precise hyperbolicity-type result for mappings of maximal rank from $\bC^2$ to very general hypersurfaces in $\bP^4$. To put our result into perspective, recall that the Kobayashi conjecture predicts the non-existence of entire curves in the general projective hypersurfaces of sufficiently high degree. As for the effective bound on the degree, for hypersurfaces in $\bP^4$ the best result to date was obtained by the second author in \cite{Rou06}, where he proved 
 the algebraic degeneracy of entire curves in a general hypersurface of degree $d \geq 593$ in $\bP^4$. Here we prove the following.
\begin{theorem}\label{thm:deg2}
Let $X_{d} \subset \bP^4$ be a very general hypersurface of degree $d \geq 93$. Then there is no holomorphic mapping of maximal rank from $\bC^2$ to $X_{d}$.
\end{theorem}

Let us remark that such a result cannot be obtained using positivity of $\Lambda^2 T_{X}^*$, since in the case of hypersurfaces $X$ of $\bP^4$ it is not hard to see that $H^0(X,\Lambda^2 T_{X}^*)=0$. Apart from the use of our constructions, there are two other ingredients in the proof of Theorem \ref{thm:deg2}. Following the recent work of Diverio, \cite{Div} we use the algebraic holomorphic Morse inequalities (see for example \cite{Tra}) to obtain non-zero global sections of $E_{2,2,m}$ vanishing on an ample divisor on $X_d$. Then    
we follow Siu's method \cite{S04}, and P\u aun's effective version of it \cite{Pau},  to produce meromorphic vector fields on jets spaces to obtain, by differentiation, other independent algebraic differential equations. These equations are enough to deduce the algebraic degeneracy of $f:\bC^2 \rightarrow X_d,$ as soon as $d\geq 93$. By Kobayashi-Ochiai the subvariety covered by the image of $f$ cannot be of general type, and then we conclude by Ein's result \cite{Ein}.

We hope that our work will allow to obtain further results in the direction of the Green-Griffiths conjecture in dimension $\geq 3$. Among the problems stemming from this work, the following seems particularly interesting to us.
\begin{pb}
Let $X$ be a smooth threefold of general type. 
Find conditions on $c_1(X),c_2(X)$ and $c_3(X)$ insuring the algebraic degeneracy of maximal rank holomorphic map from $\C^2$ to $X$.
\end{pb}

The plan of the paper is as follows. The constructions of the jet bundles $X_{k}^p$ and of the jet differentials bundles $E_{p,k,m}$ is presented in \S 2, while the natural degeneracy result Theorem \ref{thm:gen} is proved in section 3. 
The jet bundles come equipped with natural 
line bundles $\cO_{X^p_k}(\bf{a})$ whose positivity
properties are investigated in \S4 (as in \cite{Div}, these positivity properties are the key to prove existence of global section of $E_{p,k,m}$). Then we show that for certain $k$ and $p$ there cannot exist non-zero sections of $E_{p,k,m}$
(this is a generalization to all $1\leq p\leq \dim(X)$ of  \cite[Theorem 1]{Div}). In section 6 we present the first two 
applications of our general constructions. First (cf. Theorem \ref{thm:bloch}) we characterize the image of maximal rank holomorphic mappings from $\C^p$ to complex (not necessarily projective) tori, generalizing \cite[Theorem 9.1]{De95}. Then in \S6.2 we study families of maximal rank holomorphic mappings into the universal degree $d$ hypersurface in $\bP^n$ or into its complement (cf. Corollaries \ref{cor:univhypers} and \ref{cor:Cunivhypers})
generalizing \cite{DPP} and  \cite[Theorem 6]{PR} to arbitrary $p$.  These results suggest a generalization of the Kobayashi conjecture, predicting the $p$-measure hyperbolicity of a general hypersurface $X_d\subset \bP^n$ of degree $d\geq 2n+1-p$, for $1\leq p\leq n-1$ (see Problem \ref{pb:genkob}).  
In proving them, we establish a degeneracy result for families (cf. Theorem \ref{thm:famdeg}) which may be of independent interest. Then we turn to the proof of our Theorem \ref{thm:deg2}, which will occupy the final section \S7.

%
\section{Jets of holomorphic maps and generalized Demailly-Semple jet bundles}\label{S:jets}
%
\subsection{Grassmannian bundles and effective locus}

Let $(X,V)$ be a directed manifold, i.e. $X$ is a complex manifold of dimension $n$ and $V$ is a subbundle of rank $r$ of $T_{X}$. Let $x_{0} \in X$ and choose coordinates $(z_{1},\dots,z_{n})$ on a neighborhood $U$ of $x_{0}$ such that $V$ can be defined on $U$ by linear equations
\begin{displaymath}
V_{z}= \{v=\sum_{1\leq i \leq n} v_{j} \frac{\partial}{\partial{z_{j}}}; v_{j}=\sum_{1\leq k \leq r} a_{jk}(z)v_{k}, j=r+1,\dots,n\},
\end{displaymath}
with $a_{jk}$ holomorphic. Let us describe some local coordinates on the grassmannian bundle $G(p,V)$ of $p$-planes in $V$ for an integer $1\leq p \leq n$. A vector $v \in V_{z}$ is determined by $(v_{1},\dots,v_{r})$. A $p$-plane is therefore determined by a $r \times p$ matrix $(v_{ij})$. For $I=\{i_{1},\dots,i_{p}\} \subset \{1,...,r\}$, if $det(v_{ij})_{i\in I,1 \leq j \leq p} \neq 0$ then the $p$-plane can be uniquely represented by a $r \times p$ matrix $(u_{ij})$ whose $I$th $p\times p$ minor is the identity matrix. The corresponding affine chart can be described by coordinates
\begin{displaymath}
(z_{1},\dots,z_{n},(u_{ij})_{i \not\in I}).
\end{displaymath}

Let $f=(f_{1},\dots,f_{n}): (\bC^{p},0) \rightarrow (X,x)$ be a germ of holomorphic map tangent to $V$ i.e. $df_{t}(T_{\bC^{p}}) \subset V_{f(t)}$ for all $t$ in a neighborhood of $0$ in $\bC^{p}$. The map $f$ is uniquely determined by its initial value and the components $(f_{1},\dots,f_{r})$ because $f$ verifies the following system of partial differential equations 
\begin{equation}
\label{edp}
\frac{\partial{f_{j}}}{\partial{t_{i}}}=\sum_{1\leq k \leq r} a_{jk}(f)\frac{\partial{f_{k}}}{\partial{t_{i}}}, j=r+1,\dots,n;  1 \leq i \leq p.
\end{equation}

In the situation where $\rank(df_{0})=p$, the map $f$ lifts to a well defined holomorphic map $\tilde{f} : (\bC^{p},0) \rightarrow G(p,V)$ by taking $\tilde{f}=(f , [\frac {\partial f}{\partial t_{1}}\wedge \dots\wedge \frac {\partial f}{\partial t_{p}}])$. In coordinates we can describe it in the following way: take coordinates $t=(t_{1},\dots,t_{p})$ on $(\bC^{p},0)$ such that $f_{i_{k}}(t)=t_{k}$ for $i_{k} \in I=\{i_{1},\dots,i_{p}\} \subset \{1,...,r\}$. Then in the affine chart described above
\begin{displaymath}
\tilde{f}=(f_{1},\dots,f_{n}; (\frac{\partial{f_{i}}}{\partial{t_{j}}})_{i \not\in I}).
\end{displaymath}

If $p \geq 2$, it should be noted that in general, contrary to the case $p=1$, not all points $w \in G(p,V)$ are obtained as $\tilde{f}(0)$ where $\tilde{f}$ is the lifting of a germ $f: (\bC^{p},0) \rightarrow (X,x)$ tangent to $V$. Indeed, the equations (\ref{edp}) correspond to a Pfaffian system locally associated to $n-r$ linearly independent $1$-forms $\omega_{i}$, $1\leq i \leq n-r$. Classically (see for example \cite{BCG}) we introduce the closed exterior differential system $\mathcal{I} \subset \bigwedge^*(\Omega^1_X)$ generated by $\{\omega_{i},d\omega_{i}\}$ inside the exterior algebra $\bigwedge^*(\Omega^1_X)$ of $X$. A $p$-plane $E\subset  T_{X,x}$ is integral if $\psi_{|E}=0$ for all $\psi \in \mathcal{I}$. 

We denote by $Z \subset G(p,V)$ the \textit{effective locus} of integral elements of $\mathcal{I}$. Clearly, all points $w=\tilde{f}(0)$, where $\tilde{f}$ is the lifting of a germ $f: (\bC^{p},0) \rightarrow (X,x)$ tangent to $V$, lie in $Z$.

\begin{remark}\label{sing}
Let $\Delta$ be the unit disc.
If the map $f:\Delta^{p} \rightarrow X$, possibly not of maximal rank everywhere, verifies $\frac {\partial f}{\partial t_{1}}\wedge\dots \wedge \frac {\partial f}{\partial t_{p}}=g\cdot u_{1}\wedge\dots \wedge u_{p}$ with $g$ not identically zero and $u_{1}\wedge\dots\wedge u_{p}\neq 0$ then we define the lifting to be $(f,[u_{1}\wedge...\wedge u_{p}])$. In particular, the lifting of a maximal rank holomorphic map is well defined outside a codimension $\geq 2$ subset.
\end{remark}

\begin{example}
Take the cone parametrically defined by $f :(t_{1},t_{2})\rightarrow (t_{1},t_{1}t_{2},t_{1}t_{2}^{2})$. Then $\frac {\partial f}{\partial t_{1}}=u_{1}=(1,t_{2},t_{2}^{2}), \frac {\partial f}{\partial t_{2}}=t_{1}u_{2}=t_{1}.(0,1,2t_{2})$. So the lifting is $\tilde{f}=(f ,[u_{1}\wedge u_{2}])$.
\end{example}

\subsection{Jet bundles}\label{jets}
Let $(X_{0}^{p},V_{0})$ be a directed manifold with $\rank V_{0}=r$ and $p$ an integer $1\leq p \leq r$. We generalize the inductive process of Demailly-Semple jet bundles appearing in \cite{De95} for $p=1$. We define the 1-jet bundle $X_{1}^{p} = G(p,V_{0})$, as the Grassmannian bundle, and $V_{1} \subset T_{X_{1}^{p}}$ as

\begin{displaymath}
V_{{1}_{(x,\Lambda)}} =\{\xi \in T_{X_{1}^{p},(x,\Lambda)} | \pi _{\ast }\xi \in \Lambda\},
\end{displaymath}
where $\pi$ is the natural projection from $X_{1}^{p}$ to $X_{0}^{p}$. We have the two exact sequences
\begin{equation}\label{eq1}
0 \rightarrow T_{X_{1}^{p}/X_{0}^{p}} \rightarrow V_{1} \rightarrow S_{1} \rightarrow 0,
\end{equation}
\begin{equation}\label{eq2}
0 \rightarrow S_{1}\otimes S_{1}^* \rightarrow \pi^*V_{0}\otimes S_{1}^* \rightarrow T_{X_{1}^{p}/X_{0}^{p}}\rightarrow 0,
\end{equation}
where $S_{1}$ is the tautological subbundle on $X_{1}^{p}$.

\bigskip
As in \cite{LS} we can caracterize the effective locus $Z^p_{1} \subset X_{1}^{p}$ in the following way. Let $V_{0}^{\bot}$ be the kernel of the map $T_{X_{0}^{p}}^{*} \rightarrow V_{0}^{*}$. From the definition of the effective locus we immediately obtain the following characterization of $Z^p_1$.

\begin{proposition}\label{eff}
The effective locus $Z^p_{1}$ is the locus where the pairing
\begin{displaymath}
\overset{2}{\wedge}V_{1} \otimes \pi^{*}d(V_{0}^{\bot}) \rightarrow \mathcal{O}_{X_{1}^{p}}
\end{displaymath}
is zero. 
\end{proposition}

\bigskip
By induction we define the $k$-jet bundle $X_{k}^{p}=G(p,V_{k-1})$, $V_{k}$, $S_{k}$ and the effective locus $Z^p_{k}$, which is characterized by the previous proposition as the locus where the pairing
\begin{equation}\label{effectiv}
\overset{2}{\wedge}V_{k} \otimes \pi^{*}d(V_{k-1}^{\bot}) \rightarrow \mathcal{O}_{X_{k}^{p}}
\end{equation}
is zero.

We have the exact sequences
\begin{eqnarray}\label{eq:exseq}
0 \rightarrow T_{X_{k}^{p}/X_{k-1}^{p}} \rightarrow V_{k} \rightarrow S_{k} \rightarrow 0,
\end{eqnarray}
\begin{eqnarray}\label{eq:exseq2}
0 \rightarrow S_{k}\otimes S_{k}^* \rightarrow \pi^*V_{k-1}\otimes S_{k}^* \rightarrow T_{X_{k}^{p}/X_{k-1}^{p}}\rightarrow 0,
\end{eqnarray}
where $S_{k}$ is the tautological sub-bundle on $X_{k}^{p}$.

From the definition of $X^p_1$ and $S_1$ and from the exact sequence (\ref{eq1}) we get 
\begin{eqnarray*}
\dim X_{1}^{p}=n+p(r-p),\\
 \rank V_{1}=p(r-p+1),
\end{eqnarray*}
and therefore we deduce the formulas
\begin{eqnarray*}
 \dim X_{k}^{p}=n+p(r-p)+...+p^{k}(r-p)\\
 = n+(r-p)p\frac{(1-p^{k})}{(1-p)}, \text{ for } p>1,\\
 =n+k(r-1), \text{ for } p=1.\\
\text{rank } V_{k}=p^{k}(r-p)+p.
\end{eqnarray*}

Recall that we have the Plucker embedding
\begin{equation*}
\rho : G(p,V_{k-1})\hookrightarrow \bP(\overset{p}{\wedge} V_{k-1}).
\end{equation*}
Therefore we have a line bundle $\mathcal{O}_{X_{k}^{p}}(1)$ which is the restriction of $\mathcal{O}_{\bP(\overset{p}{\wedge}V_{k})}(1)$. It is also equal to $(\overset{p}{\bigwedge}{S_{k}})^{-1}$.

By definition there is a canonical injection $S_{k} \hookrightarrow \pi_{k-1,k}^{*}V_{k-1}$ where $\pi_{k-1,k} : X_{k}^{p} \rightarrow X_{k-1}^{p}$ is the natural projection. Therefore for $k\geq 2$ we obtain, by composition, a line bundle morphism 
\begin{displaymath}
\overset{p}{\bigwedge} S_{k} \hookrightarrow \overset{p}{\bigwedge} \pi_{k-1,k}^{*}V_{k-1} \rightarrow \overset{p}{\bigwedge} \pi_{k-1,k}^{*}S_{k-1}.
\end{displaymath}
We denote its zero divisor $D_{k} \subset X_{k}^{p}$. It corresponds to the set of $p$-planes in $V_{k-1}$ which meet $T_{X_{k-1}^{p}/X_{k-2}^{p}}$ non trivially.
We have
\begin{equation}\label{inj}
\mathcal{O}_{X_{k}^{p}}(1)= \pi_{k}^{*}\mathcal{O}_{X_{k-1}^{p}}(1)\otimes \mathcal{O}(D_{k}).
\end{equation}
Then we define 
\begin{displaymath}
X_{k}^{p,reg}=\underset{2 \leq j \leq k} \bigcap \pi_{j,k}^{-1}(X_{j}^{p} \setminus D_{j}),
\end{displaymath}
\begin{displaymath}
X_{k}^{p,sing}=\underset{2 \leq j \leq k} \bigcup \pi_{j,k}^{-1}(D_{j})=X_{k}^{p} \setminus  X_{k}^{p,reg},
\end{displaymath}
where $\pi_{j,k} : X_{k}^{p} \rightarrow X_{j}^{p}$ is the natural projection.

\subsection{Jet differentials}
Let $J_{k,p}=J_{k}(\bC^{p},X)$ be the space of $k$-jets of germs of holomorphic maps $ f : (\bC^{p},0) \rightarrow X $. We define the vector bundle of jet differentials of order $k$ and degree $m$, $E_{p,k,m}^{GG} \rightarrow X$ to be the vector bundle whose fiber are complex valued polynomials $Q(f^{\prime },f^{\prime \prime },\dots,f^{(k)})$ on the fibres of  $J_{k,p}$ of weight $(m,m,\dots,m)$ with respect to the action of $(\bC^{*})^{p}$:
\begin{equation*}
Q((\lambda_{1},\dots,\lambda_{p}).f)=\lambda_{1}^{m}\dots\lambda_{p}^{m}Q(f),
\end{equation*}
where 
\begin{eqnarray*}
(\lambda_{1},\dots,\lambda_{p}).f^{(i)}=(\lambda_{1},\dots,\lambda_{p}).(\frac{\partial^{i}f}{\partial t_{1}^{i}},\dots,\frac{\partial^{i}f}{\partial t_{1}^{i_{1}}\dots\partial t_{p}^{i_{p}}},\dots,\frac{\partial^{i}f}{\partial t_{p}^{i}})\\
=(\lambda_{1}^{i}\frac{\partial^{i}f}{\partial t_{1}^{i}},\dots,\lambda_{1}^{i_{1}}\dots\lambda_{p}^{i_{p}}\frac{\partial^{i}f}{\partial t_{1}^{i_{1}}\dots\partial t_{p}^{i_{p}}},\dots,\lambda_{p}^{i}\frac{\partial^{i}f}{\partial t_{p}^{i}}).
\end{eqnarray*}
\begin{example}
$E_{p,1,m}^{GG}=S^{m}T_{X}^{\ast }\otimes \dots\otimes S^{m}T_{X}^{\ast }$, p times.
\end{example}

Now we define the subbundle $E_{p,k,m} \subset E_{p,k,m}^{GG}$ by the condition
\begin{equation}
\label{inv}
Q(f\circ \phi)=(J_{\phi})^{m}Q(f)\circ \phi,
\end{equation}
for every $\phi \in G_{p,k}$ the group of $k$-jets of germs of biholomorphisms of $(\bC^{p},0)$, where $J_{\phi}$ denotes the jacobian of $\phi$.

\begin{example}[Jet differentials of order 1]\label{expl1}
Let us describe  $E_{p,1,m}$. We are looking for polynomials of $S^{m}T_{X}^{\ast }\otimes...\otimes S^{m}T_{X}^{\ast }$ ($p$ times) invariant under $G_{p,1}=GL_{p}$. From \cite{FH} p.228, we deduce that
\begin{equation*}
E_{p,1,m}=\Gamma^{(m,\dots,m,0,\dots,0)}T_{X}^{\ast},
\end{equation*}
with $p$ times $m$, $r-p$ times $0$ and $\Gamma$ which denotes the Schur functor.
\end{example}
\begin{rmk}\label{rmk:filtration}
Now let us describe how one can obtain a filtration on $E_{p,k,m}.$ We introduce the following notations. For $l\in \bN$, let $I_{l}=\{i=(i_{1},\dots,i_{p}) : i_{1}+\dots +i_{p}=l\}$. We choose orders on the $I_{l}$. Letting $i \in I_{k}$ and considering the expression of highest degree in the $\frac{\partial^{i}f_{j}}{\partial t_{1}^{i_{1}}...\partial t_{p}^{i_{p}}}$ that occurs, we see that there is a filtration
$$ F_{i}=S_{0} \subset S_{1} \dots \subset S_{[m/k]}=E_{p,k,m},$$
where
$$S_{q}/S_{q-1}\cong S^{q}T_{X}^{*}\otimes G_{i,q},$$
with $F_{i}$ which consists of polynomials of $E_{p,k,m}$ which do not involve $\frac{\partial^{i}f}{\partial t_{1}^{i_{1}}...\partial t_{p}^{i_{p}}}$ and $G_{i,q}$ consisting of polynomials of weight $(m,m,\dots,m)-qi$ which do not involve $\frac{\partial^{i}f}{\partial t_{1}^{i_{1}}...\partial t_{p}^{i_{p}}}.$ Inductively we get a filtration such that the graded terms are
$$\bigotimes_{i \in I_{1}}S^{q_{i}^{1}}T_{X}^{*} \dots  \bigotimes_{i \in I_{k}}S^{q^{l}_{i}}T_{X}^{*},$$
with
$$\sum_{l=1}^{k}\sum_{i \in I_{l}}q^{l}_{i}i=(m,\dots, m).$$
\end{rmk}

\subsection{Jet differentials and jet bundles}
In this section we will take $V=T_{X}$ in the previous constructions. Let $\Delta$ be the unit disk and $f : \Delta ^{p} \rightarrow X$ be a holomorphic map. We say it defines a regular $k$-jet if $rank(df_{t})=p$ for every $t$. In this situation, as we have seen, $f$ lifts to a well defined holomorphic map $f_{[1]} : \Delta ^{p} \rightarrow X^p_{1}$ by taking $f_{[1]}=(f , [\frac {\partial f}{\partial t_{1}}\wedge \dots \wedge \frac {\partial f}{\partial t_{p}}])$. Inductively, it lifts to a well defined map $f_{[k]} :  \Delta ^{p} \rightarrow X^p_{k}$. Moreover $\overset{p}{\wedge} df_{[k-1]} $ gives rise to a section $\overset{p}{\wedge}T_{\Delta ^{p}} \rightarrow f_{[k]}^{*}\mathcal{O}_{X^p_{k}}(-1)$.

\begin{theorem} Let $X$ be a complex manifold, $\pi_{0,k} : X^p_{k} \rightarrow X$ 
the natural projection,   $J_{k,p}^{reg} \subset J_{k,p}$ the space of regular $k$-jets and $G_{p,k}$ the group of $k$-jets biholomorphisms of $(\C^p,0)$. There is a holomorphic embedding 
\begin{displaymath}
J_{k,p}^{reg}/G_{p,k} \hookrightarrow Z^p_{k} \subset X_{k},
\end{displaymath}
which identifies $J_{k,p}^{reg}/G_{p,k}$ with $Z^{p,reg}_{k}=Z_{k}\cap X_{k}^{reg}$.
So the effective locus $Z^p_{k}$ is a natural relative compactification of $J_{k,p}^{reg}/G_{p,k}$ over $X$.
\end{theorem}

\begin{proof}
As explained at the beginning of this section, for $f \in J_{k,p}^{reg}$ we have a well defined lifting $f_{[k]}$ in $X_{k}$ and the lifting process commutes with the reparametrization i.e. $(f\circ \phi )_{[k]}=f_{[k]} \circ\varphi $. Therefore we have a map $J_{k,p}V^{reg}/G_{p,k} \rightarrow  Z^p_{k}$ defined by $ f \mod  G_{p,k} \mapsto f_{[k]}(0)$.
We can describe it in coordinates. Let $(z_{1},\dots,z_{n})$ be coordinates near a point $x \in X$ such that $V_{x}={\textrm {Vect}}(\frac{\partial }{\partial z_{1}},\dots,\frac{\partial }{\partial z_{r}})$. There exist $\{i_{1},\dots,i_{p}\} \subset \{1,\dots,r\}$ such that $g=(f_{i_{1}},\dots,f_{i_{p}})$ is a biholomorphism at $0\in \bC^{p}$. Then we take the reparametrization $\varphi=g^{-1}$. Let us suppose $\{i_{1},\dots,i_{p}\}=\{1,2,\dots,p\}$. Therefore $f\circ\varphi=(t_{1},\dots,t_{p},h_{p+1},\dots,h_{n})$. Let $h=(h_{p+1},\dots,h_{r})$. Recall that $X^p_{k}$ is a tower of $G(p, p^{i}(r-p)+p)$-bundles ($0\leq i \leq k-1$), and on the affine charts of these Grassmannians,  $f_{[k]}(0)$ is given by the collection of partial derivatives
\begin{displaymath}
((\frac{\partial h}{\partial t_{1}},\dots,\frac{\partial h}{\partial t_{p}})(0);(\frac{\partial ^{2}h}{\partial t_{i}\partial t_{j}})(0);\dots;(\frac{\partial ^{k}h}{\partial t_{i_{1}}\dots\partial t_{i_{k}}})(0)).
\end{displaymath}
This proves the holomorphy. 
As for the identification of  $J_{k,p}^{reg}/G_{p,k}$ with $Z^{p,reg}_{k}=Z^p_{k}\cap X^{p,reg}_{k}$, we refer to \cite{ASS} and \cite{LS} where it is proved that one can define the effective locus as  the closure of all possible $k$-th lifts $f_{[k]}$ of regular jets.
\end{proof}

\begin{theorem}\label{thm:imd}
Suppose $\binom{n}{p}\geq 2$ where $n=\dim X$. Then there is an injection of sheaves
\begin{displaymath}
(\pi_{0,k})_{*} \mathcal{O}_{Z^p_{k}}(m)  \subset \mathcal{O}(E_{p,k,m}).
\end{displaymath}
\end{theorem}
\begin{proof}
Let $\sigma$ be a section of $H^0(Z^p_{k,x},\mathcal{O}_{Z^p_{k,x}}(m))$ over a fiber $Z^p_{k,x}$, for $x \in X$. Let $f \in J_{k,p,x}^{reg}$, so $\overset{p}{\wedge} df_{[k-1]}(0) $ is an element of $\mathcal{O}_{X^p_{k}}(-1)$ at $f_{[k]}(0)$. Then we define
\begin{equation}\label{form1}
Q_\sigma(f',f'',\dots,f^{k})=\sigma(f_{[k]}(0)).(\overset{p}{\wedge} df_{[k-1]}(0))^{m}. 
\end{equation}
$Q_\sigma$ is holomorphic on $J_{k,p,x}^{reg}$ since $\sigma$ is. If we reparametrize $f$ by $\phi \in G_{p,k}$, then $f_{[k]}(0)$ does not change and $d(f\circ \phi)_{[k-1]}(0)=df_{[k-1]}(0).d\phi(0)$. Therefore 
\begin{displaymath}
\overset{p}{\wedge} d(f\circ \phi)_{[k-1]}(0)=J_{\phi} .\overset{p}{\wedge} df_{[k-1]}(0)
\end{displaymath}
and the invariance condition (\ref{inv}) is satisfied. Since $J_{k,p}^{reg}$ is the complement of a codimension $\binom{n}{p}$ subvariety, by the hypothesis $\binom{n}{p}\geq 2$ we have that $Q_\sigma$ extends holomorphically to $J_{k,p}$. So $Q_\sigma$ can be written as a power series in the partial derivatives of $f$ and the invariance condition implies that $Q_\sigma$ is a polynomial of weight $(m,m,\dots,m)$ for the action of $(\bC^{*})^{p}$. So there is a natural map from $(\pi_{0,k})_{*} \mathcal{O}_{Z^p_{k}}(m)$ to $\mathcal{O}(E_{p,k,m})$. To see that it is injective it is sufficient to remark that the formula (\ref{form1}) implies that it is injective over $Z^{p,reg}_{k}$.
\end{proof}

A natural problem is then to ask if the map of the previous theorem is an isomorphism (this is known to be true for $p=1$ see \cite{De95}, but we were unable to generalize this result to arbitrary $p$).
\begin{pb}
Is there an isomorphism $(\pi_{0,k})_{*} \mathcal{O}_{Z^p_{k}}(m)  \simeq \mathcal{O}(E_{p,k,m})$ ?
\end{pb}

%
\section{Algebraic degeneracy of holomorphic maps}\label{S:deg}
%

We start by recalling without proof the $p$-dimensional Ahlfors-Schwarz lemma.
\begin{lemma}[see \cite{De95}, Lemma 3.2]\label{lem:AS}
Let $\gamma(\t)=i\Sigma\gamma_{jk}dt_j\wedge d\overline t_k$ be an almost everywhere positive hermitian form on the ball $B(\0, R)\subset \C^p$ of radius $R$ such that 
$$
 i\partial\overline\partial \log \det \gamma 
 \geq \a\cdot \gamma
$$
for some constant $\a>0.$ Then 
$$
 \det (\gamma) \leq \Big( \frac{p+1}{\a R^2}\Big)^p
 \frac{1}{(1-\frac{|\t|^2}{R^2})^{p+1}}.
$$
\end{lemma}

Using the previous lemma we can prove the natural degeneracy result, which generalizes \cite[Theorem 7.8]{De95}. Before stating the theorem we recall the following standard definitions.

\begin{definition}  Let $M$ be a complex manifold and $L$ a holomorphic line bundle on it. 
\begin{enumerate}
\item[(i)]
A singular metric $h$ on $L$ is a metric which on any local trivialization 
$L_{|U}\cong U\times \C$ of $L$ is of the form 
$$
 |\xi|^2_h =|\xi|^2 \exp(-\phi),
$$
for some real valued weight function $\phi\in L^1_{loc}(U)$.
\item[(ii)] The curvature current of $L$ w.r.t. a singular  metric $h$ is the closed $(1,1)$
 current 
 $$ 
  \Theta_h(L):= \frac{i}{2\pi} \partial\overline\partial \phi
 $$
 computed in the sense of distributions.
\item[(iii)] The singularity set $\Sing(h)$ of a singular metric $h$ on $L$ is the closed subset of points in a neighborhood of which the weight $\phi$ is not locally bounded.
\item[(iv)]  A singular metric $h$ on $L$ has positive 
curvature if there exists a smooth positive definite hermitian metric $\omega$ and a continuous positive function $\eps$ on $X$ such that
$$
 \Theta_h(L)\geq \eps \cdot \omega.
$$ 
\end{enumerate}
\end{definition}

\begin{thm}\label{thm:deg}
Let $(X,V)$ be a compact complex directed manifold, and $p\geq 1$ an integer. If there exists
a metric $h=h_{p,k}$ with negative curvature on the line bundle $\cO_{X^p_k}(-1)$ on $X^p_k$, then every maximal rank holomorphic mapping $f:\C^p\lra X$ which is tangent to $V$ 
is such that  
$$
\overline{ \Im (f_{[k]})}  \subset \Sing(h),
$$
where $\overline{ \Im (f_{[k]})}$ denotes the Zariski closure of the image of $f_{[k]}$.  
\end{thm} 
\begin{proof} By the negativity of the metric $h$, there exists a smooth hermitian metric $\omega=\omega_k$ on $X_k$ and a positive function $\eps$
such that
$$
 \Theta_{h^{-1}}( \cO_{X^p_k}(1)) (\xi)\geq \eps \cdot |\xi|^2_\omega,\ \ \ \forall\xi \in V_k.
$$
On the other hand, since by (\ref{eq:exseq}) the differential 
$(\pi_k)_*$ maps $V_k$ continuously onto $\cO_{X^p_k}(-1)$, there exists a constant $C>0$ such that
$$
 |(\pi_k)_*\xi|^2_h\leq C\cdot |\xi|^2_\omega, \ \ \ \forall\xi \in V_k.
$$ 
Combining the previous two inequalities we find 
\begin{equation}\label{eq:star}
\Theta_{h^{-1}}( \cO_{X^p_k}(1)) (\xi)\geq \frac{\eps}{C} \cdot |(\pi_k)_*\xi|^2_h , \ \ \ \forall\xi \in V_k.
\end{equation} 
By abuse of notation we still denote by $f:B(\0,R)\lra X$ the restriction to the ball
of radius $R$ of a maximal rank holomorphic mapping $f:\C^p\lra X$. The Jacobian of the 
$(k-1)$-st lifting of $f$ yields a section 
$$
 F:=\Jac(f_{[k-1]}) : T_{B(\0,R)}\lra f^*_{[k]} \cO_{X^p_k}(-1)
$$
which we use to pull-back the metric $h$ to $B(\0,R)$ :
$$
 \gamma := \gamma_0(\t) \cdot d\t \otimes d\overline\t = F^*h.
$$
Notice that the metric $\gamma$ is everywhere defined for the following reason: by Remark \ref{sing} the map $f_{[k-1]}$ is well defined on a codimension $\geq 2$ subset.  Therefore $F$ is well defined on the complement of a codimension $\geq 2$ subset and we can extend it to the whole ball $B(\0,R)$.
We now argue by contradiction and assume that $f_{[k]} (B(\0,R)) \not\subset \Sing(h).$ Outside $f_{[k]}^{-1}(\Sing(h))$ using (\ref{eq:star}) we obtain: 
$$
 \frac{i\partial\overline\partial\log\det \gamma(\t)}{\gamma(\t)} = \frac{-2\pi f_{[k]}^* \Theta_h(\cO_{X^p_k} (-1))}{F^*h} =  \frac{\Theta_{h^{-1}}( \cO_{X^p_k}(1)) (\Jac(f_{[k]}(\t)))}{|\Jac(f_{[k-1](\t)})|^2_h}\geq \frac{\eps}{C},
$$
since $\Jac(f_{[k-1]})(\t)=(\pi_k)_* \Jac (f_{[k]})(\t)$.
We can therefore invoke Lemma \ref{lem:AS} and deduce the inequality 
$$
 |\Jac(f_{[k-1]}(\t))|^2_h\leq \frac{2C}{\eps}\frac{R^{-2}}{(1-|\t|^2/R^2)^2}.
$$
If we now let $R\lra+\i$, we obtain that $f_{[k-1]}$,
and therefore $f$, cannot have maximal rank, which contradicts the hypothesis.
\end{proof}
  
From  Theorem \ref{thm:deg} we deduce the following.
\begin{cor}\label{cor:deg}
Let $(X,V)$ be a compact complex directed manifold, and $p\geq 1$ an integer. Assume there exist integers $k,m>0$ and an ample line bundle $A$ on $X$ such that the linear system 
$|\cO_{X^p_k}(m)\otimes \pi_{0,k}^*\cO_X(-A)|$
is not empty. Then for every maximal rank holomorphic mapping $f:\C^p\lra X$ which is tangent to $V$ 
we have 
$$
 \overline{\Im f_{[k]}}  \subset \Base \Big(|\cO_{X^p_k}(m)\otimes \pi_{0,k}^*\cO_X(-A)|\Big).
$$
\end{cor}
\begin{proof}
Recall that a base of global sections of $\cO_{X^p_k}(m)\otimes \pi_{0,k}^*\cO_X(-A)$ allows to construct a singular metric $h$ on $\cO_{X^p_k}(-1)$ with negative curvature and such that 
$\Sing(h)=\Base(|\cO_{X^p_k}(m)\otimes \pi_{0,k}^*\cO_X(-A)|) \cup X^{\Sing}_k.$ (see \cite[\S 7.7]{De95} for the details).
Thus from Theorem \ref{thm:deg} we deduce the inclusion 
$$
\overline{ \Im (f_{[k]})} \subset \Base \Big(|\cO_{X^p_k}(m)\otimes \pi_{0,k}^*\cO_X(-A)|\Big)\cup X^{\Sing}_k.
$$
Since $f$ has maximal rank, the (Zariski closure of the) image of its $k$-th 
lifting cannot be entirely contained in $X^{\Sing}_k$,
and we are done.
\end{proof}

From the previous corollary we see that an important step towards the algebraic degeneracy of holomorphic maps is the existence of non-zero global sections in $$H^0(Z^p_{k},\cO_{Z^p_k}(m)\otimes \pi_{0,k}^*\cO_X(-A))\subset H^0(X,E_{p,k,m}\otimes \cO_X(-A)).$$
\begin{example}
The first step towards the existence of non-zero global differential operators often consists in a Riemann-Roch computation. Following example \ref{expl1}, if $X$ is a complex compact manifold of dimension $3$, a Riemann-Roch computation gives:
\begin{equation*}
\chi(X,E_{p,1,m})=\chi(X,\Gamma^{(m,m,0)}T_{X}^{\ast}) \sim \frac{1}{120}(4c_1^{3}-3c_{1}c_{2}-c_{3})m^{5}.
\end{equation*}
Since in general it may be difficult to estimate the $h^2(X,E_{p,1,m})$, from the previous computation one cannot immediately deduce the existence of a non-zero section of $E_{p,1,m}$ when $\chi>0$.
Moreover, if $X\subset \bP^4$ is a smooth hypersurface we will see in \S 5 that 
$h^0(X,E_{p,1,m})=h^0(X,\Gamma^{(m,m,0)}T_{X}^{\ast})=0$.
\end{example}

We end this section recalling the following result which will be used in the application presented in \S 6.1.
\begin{prop}[\cite{De95}, Proposition 7.2, ii)]\label{prop:posdim}
Let $L$ be a holomorphic line bundle on a compact complex manifold. Let $B_m$ be the base locus of the linear system $|mL|$ and $\Phi_m=\Phi_{|mL|}$ be the meromorphic map associated to  $|mL|$. Let $\Sigma_m$ be the closed analytic subset equal to the union of $B_m$ and the set of points $x\in (X\setminus B_m)$ such that the fiber $\Phi_m^{-1}(\Phi_m(x))$ is positive dimensional. Then, if $\Sigma_m\not=X$ and $G$ is any line bundle on $X$, the base locus of $kL-G$ is contained in $\Sigma_m$ for $k\gg0$. In particular $L$ admits a singular hermitian metric $h$ with degeneration set $\Sigma_m$ and with $\Theta_h(L)$ positive definite on $X$.
\end{prop}
%
\section{Relative and absolute positivity of the line bundles $\cO_{X^p_{k}}(\bf{a})$}\label{S:pos}
%
Let $X$ be a complex projective manifold. 
For any positive integers $p$ and $k$ consider the 
jet bundles $X^p_k$ constructed in \S\ref{S:jets}, together with their natural projections 
$$
 \pi_{j,k}:= \pi_{j,k}^p: X^p_k\lra X^p_j, \ \ 0\leq j\leq k
$$ 
(the map $\pi_{k,k}$ is the identity map).
For any ${\bf a}=(a_1,\ldots, a_k)\in \Z^k$ we define a
line bundle $\cO_{X^p_k}({\bf a})$ on $X^p_k$ as follows :
$$
 \cO_{X^p_k}({\bf a}):= \bigotimes_{j=1}^{k} 
 \pi_{j,k}^*\cO_{X^p_j}({a}_j).
$$
From equation (\ref{inj}), we deduce the equality
$$
 \pi_{j,k}^*\cO_{X^p_j}(1)=\cO_{X^p_k}(1)\otimes \cO_{X^p_k}(- \pi_{j+1,k}^*D_{j+1}-\dots - D_k).
$$
Thus by putting $D_j^*= \pi_{j+1,k}^*D_{j+1}$ for $j=1,\dots,k-1$ and $D_k^*=0$, we have
$$
 \cO_{X^p_k}({\bf a})= \cO_{X^p_k}(b_k)\otimes  \cO_{X^p_k}(-{\bf b}\cdot D^*),
$$
where
$$
{\bf b}=(b_1,\dots,b_k)\in \Z^k, b_j=a_1+\dots+a_j,
$$
$$
{\bf b}\cdot D^*=\sum_{j=1}^{k-1}b_j\pi_{j+1,k}^*D_{j+1}.
$$
In particular, if ${\bf b} \in \N^k$, we get a nontrivial morphism
\begin{equation}\label{injection}
 \cO_{X^p_k}({\bf a})= \cO_{X^p_k}(b_k)\otimes  \cO_{X^p_k}(-{\bf b}\cdot D^*) \rightarrow \cO_{X^p_k}(b_k).
\end{equation}

For the applications it is very useful to determine for which ${\bf a}=(a_1,\ldots, a_k)\in \Z^k$ the line bundle $ \cO_{X^p_k}({\bf a})$ is (relatively) positive. 
In the rest of the section we omit the index $p$ to lighten the notation. 
\begin{lemma}\label{lem:nef1}
The vector bundle $S_k\otimes \cO_{X_k}(1)$ is relatively nef over $X$.
\end{lemma}
\begin{proof}
Recall that if we have a short exact sequence of vector bundles 
$$
 0\lra F\lra E\lra G\lra 0, \ \ s=\rank(F)
$$
the nefness of $\wedge^{s-1}E^*$ implies the nefness of its quotient $F\otimes(\det F)^{-1} \cong \wedge^{s-1}F^*$.  The lemma follows from this fact applied to the tautological exact sequence on the Grassmannian.
\end{proof}
\begin{prop}\label{prop:nef}
Let $X$ be a complex projective manifold 
and $p$ and $k$ two positive integers.   
Let $A_k=A^p_k$ and $B_k=B^p_k$ be the line bundles 
on $X_k=X^p_k$ defined recursively as follows:
$$
A_0:=\cO_X(2p), \ \ A_k:= \cO_{X_k}(p^2+1)\otimes \pi^*_{k-1,k} (A_{k-1}^{\otimes(p^2+2)}),
$$
and
$$
B_0:=\cO_X, \ \ B_k:= \cO_{X_k}(p^2+1)\otimes \pi^*_{k-1,k} (B_{k-1}^{\otimes(p^2+2)}).
$$
Then the line bundle
$$
 \cO_{X_k}(1)\otimes \pi^*_{k-1,k} (A_{k-1} )
$$
(respectively 
$$
 \cO_{X_k}(1)\otimes \pi^*_{k-1,k} (B_{k-1}))
$$
is nef (resp. relatively nef over $X$.)
\end{prop}
\begin{proof} In both cases the proof is by induction on $k$ and follows the same path. Therefore we present it in a unified way. Let $L_k$ be the line bundle for which we want to prove that $
 \cO_{X_k}(1)\otimes \pi^*_{k-1,k} (L_{k-1})$ verifies the desired positivity property $(P)$ (which  is the nefness in case $L_k=A_k$, and relative nefness over $X$ in case $L_k=B_k$). By definition $
 \cO_{X_1}(1)\otimes \pi^*_{0,1} (L_{0})$ is relatively nef, when $L_0=B_0$. In case $L_0=A_0$, remark that $\Omega^1_X(2)$ is nef, hence $\Omega^p_X(2p)$ is. Using the embedding $G(p,T_X)\hookrightarrow \bP(\Omega^p_X)$ we deduce
 that $\cO_{X_1}(1)\otimes \pi_{0,1}^*\cO_{X}(2p)$
 is nef. The base of the induction is now completed.
  
Dualizing the short exact sequence  (\ref{eq:exseq}) and tensoring with $L_k$ we get
\begin{eqnarray}\label{eq:dual}
0\lra S_k^*\otimes L_k\lra V_k^*\otimes L_k
\lra \Omega^1_{X_k|X_{k-1}}\otimes L_k\lra 0.
\end{eqnarray}
If $L_k$ is such that 
\begin{eqnarray}\label{eq:a&b}
(a)\ S_k^*\otimes L_k \textrm { verifies property $(P)$ and}
 \ (b)\ \Omega^1_{X_k|X_{k-1}}\otimes L_k \textrm{ verifies property $(P)$,} 
\end{eqnarray}
then, by (\ref{eq:dual}), the bundle $V_k^*\otimes L_k$ verifies property $(P)$ over $X$ (which, by definition, is equivalent to property $(P)$ for $\cO_{X_{k+1}}(1)\otimes \pi_{k,k+1}^*L_k$).

As for (a) recall that we have the injection $S_k\hookrightarrow \pi_{k,k-1}^*V_{k-1}$, or equivalently the surjection 
$$
 \pi_{k,k-1}^*V_{k-1}^*\twoheadrightarrow S^*_k.
$$
Therefore if $\pi_{k,k-1}^*V_{k-1}^*\otimes L_k$ verifies property $(P)$, then so does $S^*_k\otimes L_k$. If we chose $L_k$ such that 
\begin{equation}\label{eq:cond(a)}
 L_k\geq  \pi_{k,k-1}^*L_{k-1},
\end{equation}
then 
$$
 \pi_{k,k-1}^*V_{k-1}^*\otimes L_k\geq \pi_{k,k-1}^*(V_{k-1}^* \otimes L_{k-1})
$$
(the relation $\geq$ is the one given by the cone of  line bundles verifying property $(P)$). 
Since $V_{k-1}^* \otimes L_{k-1}$ verifies property $(P)$ by induction we are done, and (\ref{eq:a&b}), item (a) is proved.

As for (\ref{eq:a&b}),  item (b), recall that if we have a short exact sequence
$$
 0\lra E\lra F\lra G\lra 0,
$$
then we have an injection 
$$
 G\otimes \det(E) \hookrightarrow \wedge^{\rank(E)+1}F.
$$
Recall moreover that $\det(S_k\otimes S_k^*)=\cO_{X_k}$.
Therefore from (\ref{eq:exseq2}) we get an injection
\begin{eqnarray}
T_{X_k|X_{k-1}}\hookrightarrow \wedge^{p^2+1}\pi_{k-1,k}^*V_{k-1}\otimes S^*_k,
\end{eqnarray}
or, equivalently, a surjection
\begin{eqnarray}\label{eq:surj}
\wedge^{p^2+1}\pi_{k-1,k}^*V_{k-1}^*\otimes S_k
\twoheadrightarrow\Omega^1_{X_k|X_{k-1}}.
\end{eqnarray}
Notice that $S_k\otimes \cO_{X_k}(1)$ is relatively nef over $X$ by Lemma \ref{lem:nef1} and $\pi^*_{k-1,k}(V_{k-1}^*\otimes L_{k-1})$ verifies property $(P)$  by the inductive hypothesis. Hence if
$L_k$ is such that
\begin{eqnarray}\label{eq:cond(b)}
L_k\geq (\pi_{k-1,k}^*L_{k-1}\otimes \cO_{X_k}(1))^{\otimes (p^2+1)},
\end{eqnarray}
then 
\begin{eqnarray}
L_k\otimes \wedge^{p^2+1}\pi_{k-1,k}^*V_{k-1}^*\otimes S_k \geq \wedge^{p^2+1}\big(\pi_{k-1,k}^*(V_{k-1}^*\otimes L_{k-1})\otimes S_k\otimes \cO_{X_k}(1)\big).
\end{eqnarray}
Thus from the surjection (\ref{eq:surj}) tensored by $L_k$ we get
$$
 L_k \otimes \wedge^{p^2+1} (\pi_{k-1,k}^*V_{k-1}^*\otimes S_k)
 \twoheadrightarrow \Omega^1_{X_k|X_{k-1}}\otimes L_k,
$$
and we deduce that $\Omega^1_{X_k|X_{k-1}}\otimes L_k$ verifies property $(P)$. So (\ref{eq:a&b}), item (b) is proved. To conclude the proof of the proposition we need to define $L_k$ so that 
conditions (\ref{eq:cond(a)}) and (\ref{eq:cond(b)}) are satisfied. For that we set 
$$
 L_k:= \cO_k(2)\otimes \pi_{k-1,k}^*L_{k-1}^{\otimes 3},
$$ 
and we are done.
\end{proof}
\begin{rmk}\label{rmk:relamp}
Notice that the proof of Proposition \ref{prop:nef} shows that $L_k$ is relatively ample over $X$, since $\cO_k(1)\otimes \pi_{k-1,k}^*L_{k-1}$ is relatively nef over $X$ and relatively ample over $X_{k-1}$.
\end{rmk}
\begin{cor}
Let $X$ be a complex projective manifold 
and $p$ and $k$ two positive integers.  Then the line bundle 
$$
 \Big(\cO_{X_k}(1)\otimes \pi^*_{k-1,k} (B_{k-1})\Big) \otimes \Big(\pi_{0,k}^* (\cO_X(\ell)) \Big)
$$
is nef, for all $\ell\geq 2p(p^2+2)^{k-1}$. In particular the line bundle $\cO_{X_k}(1)\otimes \pi^*_{k-1,k} (B_{k-1})$
may be written as the difference of two nef line line bundles. Namely we have :
$$
 \cO_{X_k}(1)\otimes \pi^*_{k-1,k} (B_{k-1})= \Big(\cO_{X_k}(1)\otimes \pi^*_{k-1,k} (A_{k-1})\Big) \otimes \Big( \pi_{0,k}^* \cO_X(2p(p^2+2)^{k-1}) \Big)^{-1}.
$$
\end{cor}
Because of its importance for our main application, we state more explicitly the case $p=k=2$ of the previous corollary. 
\begin{cor}\label{cor:appli}
Let $X$ be a complex projective manifold.   The line bundle $\cO_{X^2_2}(1)\otimes \pi^*_{1,2} (B_{1})$
may be written as the difference of two nef line line bundles as follows :
$$
 \cO_{X^2_2}(1)\otimes \pi^*_{1,2} (B_{1})= \Big(\cO_{X^2_2}(5,1)\otimes \pi_{0,2}^* \cO_X(24) \Big) \otimes \Big( \pi_{0,2}^* \cO_X(24) \Big)^{-1}.
$$
\end{cor}

%
\section{A vanishing theorem for global sections of $E_{p,k,m}$}\label{S:nonexistence}
%
In this section we will prove the following theorem
\begin{theorem}\label{van}
Let $X \subset \bP^N$ be a smooth complete intersection. Then
$$H^0(X,E_{p,k,m}^{GG})=0,$$
for all $m \geq 1$ and $k$ such that $$\binom {k+p} {p} -1<\dim(X)/\codim(X).$$
\end{theorem}

\begin{remark}
This result generalizes  \cite[Theorem 1]{Div} where the result was proved for $p=1$.
\end{remark}

The main ingredient of the proof is the following theorem of Br\"uckmann and Rackwitz.
\begin{theorem}[\cite{BR}, Corollary at the top of p. 634]\label{van2}
Let $X \subset \bP^N$ be a smooth complete intersection of dimension $n$. Let $T$ be any Young tableau and $t_{i}$ the number of cells inside the $i$-th column of $T$. Set $$t:=\sum_{i=1}^{N-n}t_{i},$$ 
where $t_{i}=0,$ if $i>length (T).$
Then if $t<n$ one has $$H^0(Y,\Gamma^{T}T_{X}^*)=0,$$
where $\Gamma$ is the Schur functor.
\end{theorem}

Now we can prove Theorem \ref{van}.
\begin{proof}
We have seen in Remark \ref{rmk:filtration} that the bundle $E_{p,k,m}^{GG}$ admits a filtration whose graded terms are
$$\bigotimes_{i \in I_{1}}S^{q_{i}^{1}}T_{X}^{*} \dots  \bigotimes_{i \in I_{k}}S^{q^{l}_{i}}T_{X}^{*},$$
such that
$$\sum_{l=1}^{k}\sum_{i \in I_{l}}q^{l}_{i}i=(m,\dots, m).$$
Let $n=\dim(X)$. Looking at the decomposition of these terms into irreducible $Gl(T_{X}^{*})$-representations, we see that they are all of type $\Gamma^{(\lambda_{1},\dots,\lambda_{n})}T_{X}^{*}$ with $\lambda_{i}=0$ for 
$$i>\sum_{l=1}^{k}|I_{l}|.$$
Moreover we have $|I_{l}|=\binom {l+p-1} {p-1}$,  therefore $\sum_{l=1}^{k}|I_{l}|=\binom{p+k} {p} -1.$
So if $T$ is the tableau associated to the partition $\lambda_{1}+\dots+\lambda_{n}$ we get  
$$\sum_{i=1}^{N-n}t_{i}\leq \sum_{i=1}^{N-n} \binom{p+k} {p} -1=(N-n)\left(\binom{p+k} {p} -1\right)<\frac{n}{N-n}(N-n)=n.$$
Applying Theorem \ref{van2} we obtain the vanishing of global sections of each graded piece and therefore of the bundle $E_{p,k,m}^{GG}$ itself.
\end{proof}

\begin{remark}
Since $E_{p,k,m}$ is a subbundle of $E_{p,k,m}^{GG}$, we obtain of course the vanishing of global sections of  $E_{p,k,m}$, too.
\end{remark}
%
\section{First consequences of the basic constructions}\label{S:firstcons}
%
%
\subsection{Holomorphic maps to complex tori}\label{SS:bloch}
%
The purpose of this subsection is to prove the following result.
\begin{thm}\label{thm:bloch}
Let $Y$ be a complex torus, $p\geq 1$ an integer and $f:\C^p\to Y$ a holomorphic map of maximal rank. Then the (analytic) Zariski closure $\overline{f(\C^p)^{Zar}}$ is a translate of a subtorus, i.e.. of the form $a+Y'$, where $a\in Y$ is a point and $Y'\subset Y$ is a subtorus.
\end{thm}
The theorem above has been proved by Demailly \cite[Theorem 9.1]{De95} in the case $p=1$. 
When $Z$ is an abelian variety the result has been obtained in \cite{Wo} (see also \cite{Hu-Yang} for a proof via Nevanlinna theory). 
Thanks to our general constructions, and in particular to the degeneracy result Theorem \ref{thm:deg}, Demailly's approach (which follows ideas going back to Green and Griffiths) may be adopted for arbitrary  $p\geq 1$. Notice that Theorem \ref{thm:bloch} is not a consequence of the case $p=1$, 
since a priori a union of translates of subtori may well be not isomorphic to a subtorus. 
\begin{proof}[Proof of Theorem \ref{thm:bloch}]
The proof is by induction on the dimension of $Y$.
The statement is true when $\dim(Y)=1$. Suppose the theorem holds in dimension $n-1$.
Let $Y$ be an $n$-dimensional complex torus.
Let $f:\C^p\to Y$ be a maximal rank holomorphic map. Let $X$ be the Zariski closure of its image.  
Let $Y_k=Y_k^p$ be the $k$-jet bundle of $Y$ and $X_k$ the closure of $X_k^{reg}=X_k^{p, reg}$ in $Y_k$. Since the tangent bundle $T_Y$ is trivial, so is the $k$-jet bundle $Y_k=Y\times \bG_k$, where $\bG_k$ is the tower of Grassmannians over an arbitrary point of $Y$ which we have  defined in \S 2.  
By Remark \ref{rmk:relamp} there exists a weight ${\bf a}\in \N^k$ for which $\cO_{Y_k}({\bf a})$ is relatively very ample, i.e.. $\cO_{Y_k}({\bf a})$ is the pull back of a very ample line bundle $\cO_{\bG_k}({\bf a})$ on $\bG_k$, via the natural projection $Y_k\times \bG_k \to \bG_k$.  Let $\Phi_k:X_k\to \bG_k$ 
be the morphism induced by the second projection. 
By functoriality we have 
$$
 \cO_{X_k}({\bf a})=\cO_{Y_k}({\bf a})_{|X_k}=\Phi_k^*\cO_{\bG_k}({\bf a}).
$$
Let $\Sigma_k\subset X_k$ be the set of points $x\in X_k$ such that the fiber $\Phi_k^{-1}(\Phi_k(x))$ is positive dimensional. Assume $\Sigma_k\not=X_k$.
Since the base locus of $\cO_{X_k}({\bf a})$ is empty, by Proposition \ref{prop:posdim} the line bundle $\cO_{X_k}({\bf a})$ carries a hermitian metric with degeneration set $\Sigma_k$ and with strictly positive definite curvature. By Theorem \ref{thm:deg} we deduce that 
$$
 f_{[k]}(\C^p) \subset \Sigma_k
$$
(this inclusion is obvious if $\Sigma_k=X_k$).
Because of the previous inclusion,  through every 
point $f_{[k]}({\bf t_0})$ of the image of $f_{[k]}$ there is a germ of a (positive dimensional) variety in the fiber $\Phi_k^{-1}(\Phi_k(f_{[k]}({\bf t_0})))$. Let ${t'}\mapsto u({ t'})=(y({ t'}), j_k)\in X_k\subset Y\times \bG_k$ be the germ of a curve in  $\Phi_k^{-1}(\Phi_k(f_{[k]}({\bf t_0})))$. We have $u({0})= f_{[k]}({\bf t_0})=(y_0,j_k)$ and $y_0=f({\bf t_0})$.
Then $(y({t'}), j_k)$ is the image of $f_{[k]}({\bf t_0})$ by the $k$-th lifting of the translation $\tau_s:y\mapsto y+s$, where $s=y({ t'})- y_0$. Notice
that the image of $f$  cannot be entirely contained in the singular locus $X^{Sing}$ of $X$. Therefore we may choose ${\bf t_0}$ so that $f({\bf t_0})$ is a regular point.  Consider
$$
A_k(f):=\{s\in Y:f_{[k]}({\bf t_0})\in X_k\cap (\tau_{-s}(X))_k  \},
$$
where $ (\tau_{-s}(X))_k$ is the $k$-th jet bundle of
the subvariety $\tau_{-s}(X)\subset Y$, obtained by the translation defined by $-s$. By construction $A_k(f)$ is an analytic subset containing the germ of the curve $t'\mapsto s(t'):=y(t')-y_0$. Since  $A_1(f)\supset A_2(f)\supset\ldots\supset  A_k(f)\supset $
by notherianity the sequence must stabilize at some $A_k(f)$. Therefore there exists a germ of curve 
$D(0,r)\to Y,\ t'\mapsto s(t')$ such that the infinite jet defined by $f$ at $\bf t_0$ is $s(t')$-translation invariant for all $t'\in D(0,r)$.  By the uniqueness of the analytic continuation the same must be true for all ${\bf t}\in \C^p$, that is
$$
 s(t')+f({\bf t})\in X,\ \forall t'\in D(0,r),\ \forall{\bf t}\in \C^p.
$$
By the irreducibility of $X$ we have
$$
 s(t')+X=X.
$$
Consider 
$$
W:=\{s\in Y: s+X=X\}.
$$
Then $W$ is a closed subgroup of $Y$ which is not empty by the above. Let $p:Y\to Y/W$ be the quotient map. Notice that $Y/W$ is a complex torus of dimension $<\dim(Y)$.  Consider the composed map $\hat f :=p\circ f:\C^p\to Y/W$. Take a subspace $\C^q\subset \C^p$ such that $\overline {\hat f (\C^q)}^{Zar}=\overline {\hat f (\C^p)}^{Zar}$
and the restriction $\hat f_{|\C^q}$ has maximal rank.
By the inductive hypothesis $\overline {\hat f (\C^q)}^{Zar}$ is the translate $\hat s+\hat T$ of a subtorus $\hat T$. Since, by definition of $W$, the variety $X$ is $W$-invariant, we get that $X=s+p^{-1}(\hat T)$.
Since $p^{-1}(\hat T)$ is a closed subgroup of $Y$ the theorem is proved.
\end{proof}
%
\subsection{Families of holomorphic maps to the universal hypersurface and to its complement}\label{SS:fam}
%
When dealing with families of varieties the following result, which is a relative version (in the case $k=1$) of the degeneracy theorem proved in \S\ref{S:deg}, turns out to be very useful.

\begin{thm}\label{thm:famdeg}
Let $T$ be a complex manifold of dimension $N$ containing $\bC^N$. 
Let $\cY\to T$ be a smooth family of projective varieties of relative dimension $n$, $\cD$ an effective divisor on $\cY$ with simple normal crossing support that 
does not contain any  fiber $Y_t$ of  the family. 
Let $\cA$ be  a relatively ample line bundle on 
$\cY$ and $m$ and $p$ positive integers. Let $U\subset T$ be an open subset and $\cY_U\to U$ be the induced family.
Then for every section
$Q\in H^0(\cY_U, \Sym^m \Omega_{\cY_U}^{p+N}(\log \cD) \otimes\cO_{\cY_U}(-\cA)),$ every holomorphic mapping $\Phi :\bC^p\times U\to \cY\setminus \cD$ ($1\leq p \leq n$) such that 
$\Phi(\bC^p\times\{t\})\subset Y_t$ for every $t\in U$ and whose Jacobian $J_\Phi$ has maximal rank at some point must satisfy the algebraic differential equation $Q(\Phi)=0.$
\end{thm}
An immediate consequence is the following.
\begin{cor}\label{cor:famdeg}
In the situation of Theorem \ref{thm:famdeg}, 
if the sheaf $\Sym^m \Omega_{\cY_U}^{p+N}(\log \cD) \otimes\cO_{\cY_U}(-\cA))$ is generated by its global sections, then there is no maximal rank
 holomorphic mapping $\Phi :\bC^p\times U\to \cY\setminus \cD$.
\end{cor}
Set $\bP^{N_d}:=\bP H^0(\cO_{\bP^n}(d))^*$. From the previous degeneracy result, using 
the global generation of the sheaf of twisted (logarithmic) vector fields on the universal  degree $d$ hypersurface $\cX\subset \bP^n\times \bP^{N_d}$ (respectively on its complement) we deduce the following corollaries.

\begin{cor}\label{cor:univhypers}
Let
$\Phi:\bC^p\times \bP^{N_d}\to {\mathcal X}$ ($1\leq p \leq n-1$) be a holomorphic map
such that $\Phi(\bC^p\times\{ t\})\subset
X_t$ for all $t\in \bP^{N_d}$. 
If $d\ge 2n +1-p$,   the rank of $\Phi$ cannot be maximal anywhere.
\end{cor}

\begin{cor}\label{cor:Cunivhypers}
Let
$\Phi:\bC^p\times \bP^{N_d}\to \bP^n\times \bP^{N_d}\setminus {\mathcal X}$ ($1\leq p \leq n$) be a holomorphic map
such that $\Phi(\bC^p\times\{ t\})\subset
\bP^n\setminus X_t$ for all $t\in \bP^{N_d}$. 
If $d\ge 2n +2-p$,   the rank of $\Phi$ cannot be maximal anywhere. 
\end{cor}
The previous corollaries were proved for $p=1$ in 
\cite[Theorem 1]{DPP}, and \cite[Theorem 6]{PR}
respectively. The compact case may be seen as a generalization to the transcendental setting of the main result of \cite{Ein}.
For the terminology and notation in the logarithmic framework we refer the reader to \cite{Ii}.
Assuming Theorem \ref{thm:famdeg} for the moment
let us show how we can derive from it the two corollaries.  To simplify the notation, line bundles of the form $pr_1^*\cO_{\bP^n}(a)\otimes pr_2^*\cO_{\bP^{N_d}}(b)$ will be denoted by $\cO(a,b)$. Recall that $\cX\subset \bP^n\times \bP^{N_d}$ is of bidegree $(d,1)$. Thus its canonical bundle is $K_{\cX}\cong \cO_{\cX}(d-n-1,-N_d)$. 
\begin{proof}[Proof of Corollary \ref{cor:univhypers}]
Notice that we have the following isomorphism
\begin{equation}\label{eq:iso}
\Omega^{N+p}_{\mathcal X}(-1,0)\cong (\wedge^{n-1-p}T_{\mathcal X})\otimes K_{\cX}(-1,0)
\end{equation}
\begin{eqnarray*}
\cong
\Big(\wedge^{n-1-p}\big(T_{\mathcal X}(1,0)\big)\Big)\otimes \cO_{\cX}(d-2n-1+p,-N_d).
\end{eqnarray*}
By \cite[Lemma 4]{S04} the sheaf $T_{\mathcal X}(1,0)$ is globally generated. Take an open subset $U\subset \bP^{N_d}$ trivializing $\cO_{\cX}(0,-N_d)$.
By (\ref{eq:iso}) the sheaf $\Omega^{N+p}_{\mathcal X_U}(-1,0)$ is generated by its global sections as soon as $d\geq 2n+1-p$. 
The corollary now follows from Corollary \ref{cor:famdeg} applied to $\cY=\cX,\ \cD=0,\ \cA=\cO_{\cX}(-1,0)$ and $T$ equal to  $\bP^{N_d}$.
\end{proof}
\begin{proof}
[Proof of Corollary \ref{cor:Cunivhypers}]
Set $\bP=\bP^n\times \bP^{N_d}.$
We have the following isomorphism
\begin{eqnarray}\label{eq:iso2}
\Omega^{N+p}_{\bP}(\log \mathcal X)(-1,0)\cong (\wedge^{n-p}T_{\bP}(-\log \mathcal X))\otimes K_{\bP}\otimes\cO_{\bP}(\cX)(-1,0)
\end{eqnarray}
\begin{eqnarray*}
\cong
\Big(\wedge^{n-p}\big(T_{\bP}(-\log \mathcal X)(1,0)\big)\Big)\otimes \cO_{\cX}(d-2n-2+p,-N_d).
\end{eqnarray*}
By \cite[Proposition 11]{PR} the sheaf $T_{\bP}(-\log \mathcal X)(1,0)$ is globally generated. 
Take an open subset $U\subset \bP^{N_d}$ trivializing $\cO_{\cX}(0,-N_d)$. By (\ref{eq:iso2}) the sheaf $\Omega^{N+p}_{\bP_U}(\log \mathcal X_U)(-1,0)$ is  generated by its global sections as soon as $d\geq 2n+2-p$. 
The corollary now follows from Corollary \ref{cor:famdeg} applied to $\cY=\bP,\ \cD=\cX,\ \cA=\cO_{\bP}(-1,0)$ and $T$ equal to $\bP^{N_d}$.
\end{proof}
We pass now to the proof of the relative degeneracy result. 
\begin{proof}
[Proof of  Theorem \ref{thm:famdeg}] 
Consider  a holomorphic map $\Phi: \bC^p\times T\to \cY\setminus \cD$ with maximal rank and a section $Q\in H^0(\cY, \Sym^m \Omega_{\cY}^{p+N}(\log \cD) \otimes\cO_\cY(-\cA))$. To lighten the (already heavy) notation we will write the proof for $m=1$. 
We let $z_1,\ldots,z_p$ be the coordinates on $\C^p$ and $\xi_1,\ldots,\xi_{N_d}$ be the coordinates on $T$.
Let 
$$
J_\Phi(\underline z, \underline \xi)=
{\frac{\partial \Phi }{\partial z_1}}\wedge {\frac{\partial \Phi
}{\partial z_2}}\wedge\dots \wedge {\frac{\partial \Phi
}{\partial z_{p}}}\wedge {\frac{\partial \Phi }{\partial \xi_1}}\wedge {\frac{\partial \Phi
}{\partial \xi_2}}\wedge\dots \wedge {\frac{\partial \Phi
}{\partial \xi_{N_d}}}(\underline z,\underline \xi)\in \bigwedge ^{p+N_d}
\Phi^*T_{\cY, \Phi(\underline z,\underline \xi)}
$$ 
be the Jacobian of $\Phi$.
We argue by contradiction and suppose that
$$
 \Phi^*Q (J_\Phi)
$$
gives a non-zero section of $\Phi^*\cO_\cY(-\cA)$.
Without loss of generality we may assume that 
$  \Phi^*Q(J_\Phi)$ is non-zero at the origin. 
We denote by $\bB(r)\subset \bC$ the disc of radius
$r$ centered at the origin.

Fix a positive $\delta_0$ and, for every positive integer $k$ consider the following
sequence of maps
$$
 \Phi_k: \bB(\delta_0k)^{N_d+ p}\to \cY \setminus{\cD}
$$
 given by
$$
\displaystyle \Phi_k(z_1,\ldots,z_p, \xi_1,\dots,\xi_{N_d})= \Phi (z_1 k^{N_d},\ldots,z_p k^{N_d},
\frac{1}{k^p}\xi_1,\dots,\frac{1}{k^p}\xi_{N_d}).
$$
Since $\cO_\cY(\cA)$ is relatively ample eventually after shrinking the base $T$ we can endow it with a metric 
$h_\cA$  with positive curvature. 
For any $\underline w\in \bB(\delta_0k)^{N_d+ p}$ set
\begin{equation}\label{f_k}
 f_k (\underline w)=
\Vert \Phi_k^*Q (J_{\Phi_k}) (\underline w) \Vert ^{2/(N_d+p)}_{\Phi_k ^*h^{-1}_\cA}.
\end{equation}
By construction we have that the value at the origin 
is independent of $k$ and not zero, since
\begin{equation}\label{eq:notzero}
   f_k(\underline{0})= \Vert \Phi_k^*Q (J_{\Phi}) (\underline 0) \Vert ^{2/(N_d+p)}_{\Phi ^*h^{-1}_\cA}\not=0.
\end{equation}
Theorem \ref{thm:famdeg} follows from the
contradiction of (\ref{eq:notzero}) with the 
following.
\begin{prop}\label{tozero}
 For each $k\geq 1$ we have $f_k(\underline{0})\leq C\cdot k^{-2}$.
 In particular, as $k\to \infty$, we have 
 $f_k(\underline{0})\to 0$.
\end{prop}
\begin{proof}
First notice that because of the positivity of $h_A$, for every $k>0$ we have
\begin{equation}\label{eq:1}
\Delta \log(f_k)\geq C\cdot f_k.
\end{equation}
Then, consider the volume form of the Poincar\'e metric on the
polydisc

$$\psi_k= \prod_{j=1}^{p}\frac{1}{\Bigl(1- \frac{\vert z_j\vert ^2}{\delta_0^2k^2}\Bigr)^2}
\prod_{j=1}^{N_d}\frac{1}{\Bigl(1- \frac{\vert \xi_j\vert
^2}{\delta_0^2k^2}\Bigr)^2}.$$ A  computation shows that
\begin{equation}\label{upperbound}
\Delta \log \psi_k\leq C\cdot k^{-2}\psi_k.
\end{equation}

Consider the function $\displaystyle (\underline z, \underline\xi)\mapsto \frac{f_k(\underline z,
\xi)}{\psi_k(\underline z, \underline\xi)}$. Its maximum cannot be achieved at a
boundary point of the domain, since $\psi_k$ goes to infinity as
$(\underline z, \xi)$ goes to the boundary. So at the maximum point $(\underline z_0,
\underline \xi_0)$, we have
$$\Delta \log f_k/\psi_k\leq 0.$$
 This inequality, combined with  (\ref{eq:1})
and (\ref{upperbound}), gives
$$f_k(\underline z_0, \underline \xi_0)\leq C\cdot k^{-2}\psi_k(\underline z_0, \underline\xi_0).$$
Since the previous inequality is verified at the maximum point of the
quotient, the same is true at an arbitrary point, thus, in
particular, at the origin:
$$
f_k(\underline{0})\leq C \cdot k^{-2}.
$$
\end{proof}
The proof of Theorem \ref{thm:famdeg} is now completed.
\end{proof}
For every $1\leq p\leq \dim(X)$, Eisenmann introduced in \cite{E} $p$-dimensional pseudo-metrics ${\bf e}^p_X$ which are generalizations of the Kobayashi pseudo-metric (see also \cite[\S 1 and 3]{De95} for definitions and basic properties of such objects). A variety is called infinitesimally $p$-measure hyperbolic if the $p$-dimensional Eisenmann pseudo-metrics ${\bf e}^p_X$ is positive definite on every fiber of $\wedge^pT_X$ (in the non-compact case one requires ${\bf e}^p_X$ to further satisfy a locally uniform lower bound in terms of any smooth metric). The major difference between the case $p=1$ and $p>1$ is that the 
relevant geometric objects in the latter case are holomophic maps $f:\bC\times D(0,R)^{p-1}\to X$.
Indeed the non-existence of maximal rank holomophic maps $f:\bC^p\to X$ does not imply the $p$-measure hyperbolicity, when $p>1.$ 
Nevertheless the previous results, together with  Theorem \ref{thm:deg2}, may be seen as first steps 
towards the $p$-measure hyperbolicity of general projective hypersurfaces of high degree (and of their complements) and
suggest the following problem, which generalizes Kobayashi's conjecture.
\begin{pb}\label{pb:genkob}
Let $X_{d}\subset \bP^n$ be a general hypersurface of degree $d$.
\begin{enumerate}
\item[(i)] If $d\geq 2n+1-p$ then $ X_{d}$ is $p$-measure hyperbolic ($1\leq p \leq n-1$). 
\item[(ii)] If $d\geq 2n+2-p$ then $\bP^n\setminus X_{d}$ is $p$-measure hyperbolic ($1\leq p \leq n$). 
\end{enumerate}

\end{pb}

%
\section{The main application: holomorphic maps from $\bC^2$ to general hypersurfaces in $\bP^4$}\label{S:P4}
%
In this section we investigate the case of holomorphic maps $f:\bC^2 \rightarrow X$ where $X$ is a hypersurface of $\bP^4$, using our construction of generalized Demailly-Semple jet bundles. This approach generalizes what was done in \cite{Rou06} for entire curves. As $p=2$, we will forget the index $p$ in the rest of the section.

As explained above, the first step is to obtain global differential operators. The first remark here is that there are no global differential operators of order $1$. Indeed, by theorem \ref{van}
$$H^0(X,E_{2,1,m})=0.$$

So we have to look for global differential operators of order $2$. This approach will make use of explicit computations on the cohomology of jet bundles.

\subsection{Cohomology of jet bundles}
First, let us recall the cohomology on Grassmannian bundles. Let $V$ be a vector bundle of rank $r$ on $X$ and $X_{1}=G(p,V)$ the Grassmannian bundle. Then $H^{*}X_{1}$ is the algebra over $H^*X$ generated by $a_{1},...,a_{p},b_{1},...,b_{r-p}$ modulo the relations
\begin{displaymath}
\sum_{0 \leq i \leq k} a_{i}b_{k-i}=c_{k}(V),
\end{displaymath}\label{eq4}
for $k=1,...,r, a_{i}=c_{i}(S), b_{j}=c_{j}(Q)$ with $S$ and $Q$ defined by
\begin{equation}
0 \rightarrow S \rightarrow \pi^{*}V \rightarrow Q\rightarrow 0,
\end{equation}
where $\pi: X_1 \rightarrow X$ is the natural projection and $S$ the tautological bundle.

We make the computations in our situation where $X$ is  a hypersurface of $\bP^4$ and $V_{0}=T_{X}$, using the notations of section \ref{jets} where we have defined $X_k$, $V_k$ and $S_k.$

\bigskip
On $X_{1}=G(2,V_{0})$, we have
\begin{equation}
0 \rightarrow S_{1} \rightarrow \pi^{*}V_{0} \rightarrow Q_{1}\rightarrow 0,
\end{equation}
with $\rank S_{1}=2, \rank V_{0}=3, \rank Q_{1}=1.$
The relations defining the cohomology are:
\begin{eqnarray*}
a_{1}+b_{1} & = & c_{1}, \\
a_{2}+a_{1}b_{1} & = & c_{2},\\
a_{2}b_{1}&=&c_{3},
\end{eqnarray*}
where $a_{i}, b_{i}$ and $c_{i}$ are the $i$-th Chern classes of $S_{1}, Q_{1}$ and $V_{0}$, respectively.

Eliminating the $b_{i}$'s we obtain the relation
\begin{equation}
\label{rel1}
a_{1}^3-2a_{1}^2c_{1}+a_{1}c_{1}^2+a_{1}c_{2}-c_{1}c_{2}+c_{3}=0.
\end{equation}

\bigskip
On $X_{2}=G(2,V_{1})$, we have
\begin{equation}
0 \rightarrow S_{2} \rightarrow \pi^{*}V_{1} \rightarrow Q_{2}\rightarrow 0,
\end{equation}
with $\rank S_{2}=2, \rank V_{1}=4, \rank Q_{2}=2.$
The relations defining the cohomology are:
\begin{eqnarray*}
d_{1}+e_{1} & = & f_{1}, \\
d_{2}+d_{1}e_{1}+e_{2} & = & f_{2},\\
d_{1}e_{2}+d_{2}e_{1}&=&f_{3},\\
d_{2}e_{2}&=&f_{4}.
\end{eqnarray*}
where $d_{i}, e_{i}$ and  $f_{i}$ are the $i$-th Chern classes of $S_{2}, Q_{2}$ and $V_{1}$ respectively.

Eliminating the $e_{i}$'s we obtain the relations
\begin{eqnarray}
\label{rel2}
-{d_{{1}}}^{6}+3\,{d_{{1}}}^{5}f_{{1}}-3\,{d_{{1}}}^{4}{f_{{1}}}^{2}-2
\,{d_{{1}}}^{4}f_{{2}}+4\,{d_{{1}}}^{3}f_{{2}}f_{{1}}
&&\nonumber\\+{d_{{1}}}^{3}{f_
{{1}}}^{3}-f_{{3}}{d_{{1}}}^{2}f_{{1}}-{f_{{2}}}^{2}{d_{{1}}}^{2}+4\,{
d_{{1}}}^{2}f_{{4}}-2\,{d_{{1}}}^{2}{f_{{1}}}^{2}f_{{2}}+f_{{3}}d_{{1}
}{f_{{1}}}^{2}
&&\\-4\,d_{{1}}f_{{4}}f_{{1}}+{f_{{2}}}^{2}d_{{1}}f_{{1}}-f_
{{3}}f_{{2}}f_{{1}}+{f_{{3}}}^{2}+f_{{4}}{f_{{1}}}^{2}&=&0,\nonumber
\end{eqnarray}

\begin{eqnarray}\label{rel3}
-{d_{{1}}}^{5}-f_{{3}}f_{{2}}+3\,f_{{4}}f_{{1}}+f_{{3}}{f_{{1}}}^{2}-2
\,d_{{2}}f_{{3}}-2\,d_{{2}}f_{{1}}f_{{2}}-d_{{2}}{f_{{1}}}^{3}
&&\nonumber\\+5\,f_{{
1}}{d_{{2}}}^{2}+4\,d_{{1}}f_{{4}}+d_{{1}}{f_{{2}}}^{2}+f_{{3}}f_{{1}}
d_{{1}}-d_{{1}}f_{{2}}{f_{{1}}}^{2}-4\,f_{{2}}d_{{2}}d_{{1}}+4\,d_{{2}
}{f_{{1}}}^{2}d_{{1}}
&&\\+f_{{3}}{d_{{1}}}^{2}-2\,{d_{{1}}}^{2}f_{{2}}f_{{
1}}+{f_{{1}}}^{3}{d_{{1}}}^{2}&=&0.\nonumber
\end{eqnarray}

Relations \ref{rel1}, \ref{rel2} and \ref{rel3} will be used in the next section.

\bigskip
Now, we compute the Chern classes of $V_{1}$. From the exact sequences (\ref{eq1}) and (\ref{eq2}), we deduce
\begin{displaymath}
ch(V_{1})=ch(S_{1})+ch(\pi^*V_{0}\otimes S_{1}^*)-ch(S_{1}\otimes S_{1}^*).
\end{displaymath}

We have
\begin{eqnarray}
ch(S_{1}) & = & 2+a_{{1}}+1/2\,{a_{{1}}}^{2}-a_{{2}}, \\
ch(\pi^*V_{0}) & = & 3+c_{1}+\frac{1}{2}(c_{1}^2-2c_{2})+\frac{1}{6}(c_{1}^3-3c_{1}c_{2}+3c_{3}).
\end{eqnarray}
The computation gives
\begin{eqnarray*}
ch(V_{1})&=&4+1/3\,{c_{{1}}}^{3}-2\,c_{{2}}-c_{{1}}a_{{1}}+1/2\,c_{{1}}{a_{{1}}}^{
2}-c_{{1}}a_{{2}}-1/2\,{c_{{1}}}^{2}a_{{1}}\\
&&+1/4\,{c_{{1}}}^{2}{a_{{1}}
}^{2}-1/2\,{c_{{1}}}^{2}a_{{2}}+c_{{2}}a_{{1}}-1/2\,c_{{2}}{a_{{1}}}^{
2}+c_{{2}}a_{{2}}-1/6\,{c_{{1}}}^{3}a_{{1}}+1/12\,{c_{{1}}}^{3}{a_{{1}
}}^{2}-1/6\,{c_{{1}}}^{3}a_{{2}}\\
&&-1/2\,c_{{3}}a_{{1}}+1/4\,c_{{3}}{a_{{
1}}}^{2}-1/2\,c_{{3}}a_{{2}}+{a_{{1}}}^{2}a_{{2}}+2\,c_{{1}}+1/2\,c_{{
1}}c_{{2}}a_{{1}}-1/4\,c_{{1}}c_{{2}}{a_{{1}}}^{2}-c_{{1}}c_{{2}}+{c_{
{1}}}^{2}\\
&&+1/2\,c_{{1}}c_{{2}}a_{{2}}-1/4\,{a_{{1}}}^{4}-2\,a_{{1}}+{a_
{{1}}}^{2}-{a_{{2}}}^{2}+c_{{3}}.
\end{eqnarray*}

Then we deduce
\begin{proposition}\label{chern} Let $X\subset \bP^4$ be a hypersurface, $V_0=T_X$, $V_1$, $S_1$ defined as in section \ref{jets}, $c_i$ the Chern classes of $X$ and $a_i$ the Chern classes of $S_1$. Then
\begin{eqnarray*}
c_{1}(V_{1})&=&2\,c_{{1}}-2\,a_{{1}},\\
c_{2}(V_{1})&=&{c_{{1}}}^{2}-3\,c_{{1}}a_{{1}}+{a_{{1}}}^{2}+2\,c_{{2}},\\
c_{3}(V_{1})&=&-{c_{{1}}}^{2}a_{{1}}+c_{{1}}{a_{{1}}}^{2}+2/3\,{a_{{1}}}^{3}+2\,c_{{1
}}c_{{2}}-2\,c_{{2}}a_{{1}}-2\,c_{{1}}a_{{2}}+2\,c_{{3}},\\
c_{4}(V_{1})&=&-c_{{3}}a_{{1}}-6\,{a_{{1}}}^{2}a_{{2}}+c_{{2}}{a_{{1}}}^{2}-{c_{{1}}}
^{2}a_{{2}}-6\,c_{{2}}a_{{2}}+4\,c_{{1}}c_{{3}}+1/3\,c_{{1}}{a_{{1}}}^
{3}\\
&&+1/2\,{c_{{1}}}^{4}-2\,{c_{{1}}}^{2}c_{{2}}+2\,{c_{{2}}}^{2}+6\,{a_
{{2}}}^{2}+2/3\,{a_{{1}}}^{4}-c_{{1}}c_{{2}}a_{{1}}+4\,a_{{1}}c_{{1}}a
_{{2}}.
\end{eqnarray*}
\end{proposition}

\subsection{Existence of global 2-jet differentials}
We follow the strategy described in \cite{Div} to prove the existence of non zero global section of $E_{2,2,m}$. The main point here is that we want to find non trivial differential equations i.e. sections non identically zero on the effective locus. We prove the following.
\begin{theorem}\label{jdiff}
Let $X\subset \bP^4$ be an hypersuface of degree $d\geq 19$ and $A$ an ample line bundle. Then there exists non-zero global sections in
\begin{displaymath}
H^0(Z_{2},\mathcal{O}_{Z_{2}}(m)\otimes \pi_{2,0}^*A^{-1}) \subset H^0(X,E_{2,2,m}\otimes A^{-1})\neq 0,
\end{displaymath} 
for $m$ large enough. Therefore, by Corollary \ref{cor:deg}, any maximal rank holomorphic map $f: \bC^2 \rightarrow X$ satisfies the corresponding algebraic differential equation.
\end{theorem}

Let $F_{2}=\cO_{X^2_2}(5,1)\otimes \pi_{0,2}^* \cO_X(24)$ and $G_{2}=\pi_{2,0}^*\cO_{X}(24)$. Recall from Corollary \ref{cor:appli} that $L_{2}=F_{2}\otimes G_{2}^{-1}$ is a difference of two nef line bundles. To prove the above theorem it is enough to prove the following.

\begin{proposition}\label{big}
Under the hypotheses of Theorem \ref{jdiff}, the line bundle $L_{2}|_{Z_{2}}$ is big.
\end{proposition}

Let $u_{i}=c_{1}(\cO_{X_{i}}(1)).$ We will need the following. 
\begin{lemma}\label{lem:Z_2}
$Z_{2}\subset X_{2}$ is a divisor whose class  in $\Pic(X_{2})$ is given by
\begin{displaymath}
Z_{2}=u_{2}+u_{1}-\pi_{2,0}^*K_{X}.
\end{displaymath}
\end{lemma}
\begin{proof}[Proof of Lemma \ref{lem:Z_2}]
Recall from relation \ref{effectiv} that the effective locus $Z_{2}$ is the locus where the pairing
\begin{displaymath}
\overset{2}{\wedge}V_{2} \otimes \pi^{*}d(V_{1}^{\bot}) \rightarrow \mathcal{O}_{X_{2}}
\end{displaymath}
is zero, where $\pi= \pi^{*}_{2,1}.$ Equivalently, it is the locus of zeros of the map
\begin{displaymath} 
\pi^{*}d(V_{1}^{\bot}) \rightarrow \overset{2}{\wedge}S_{2}.
\end{displaymath}
Since $V_{1}^{\bot}$ is the kernel of $T_{X_{1}}^* \rightarrow V_{1}^*$, it is of rank $1$. Locally, it is generated by a $1$-form $\omega$.
As an $\mathcal{O}_{X_{1}}$-module $d(V_{1}^{\bot})$ is generated by $d\omega$ and the forms $df\wedge \omega$. Consider $L= d(V_{1}^{\bot})/ (T_{X_{1}}^*\wedge V_{1}^{\bot})$. Then we see that $Z_{2}$ is the locus where the pairing
\begin{displaymath}
\overset{2}{\wedge}V_{2} \otimes \pi^{*} L \rightarrow \mathcal{O}_{X_{2}}
\end{displaymath}
is zero. Indeed, take $v_{1},v_{2} \in V_{2}$ then 
\begin{displaymath}
<\pi^{*} (df\wedge \omega), v_{1}\wedge v_{2}>=<df\wedge \omega, \pi_{*} (v_{1}\wedge v_{2})>=0,
\end{displaymath}
since $\pi_{*} v_{i}\in V_{1} $ and therefore $\omega(\pi_{*} v_{i})=0$. 

Remark that $L$ is a locally free module of rank $1$, locally generated by $d\omega$, and in $Pic(X_{1})$, $L=V_{1}^{\bot}$. Indeed, $V_{1}^{\bot}$ is given by the transition functions $g_{ij}$ such that $\omega_{i}=g_{ij}\omega_{j}$ where $\omega_{i}$ defines locally   $V_{1}^{\bot}$ on $U_{i}$. Now, the transition functions of $L$ are given by $\overline{d\omega_{i}}=\overline{g_{ij}d\omega_{j}+dg_{ij}\wedge \omega}=g_{ij}\overline{d\omega_{j}}$ in $d(V_{1}^{\bot})/ (T_{X_{1}}^*\wedge V_{1}^{\bot})$.

Now, we can characterize $Z_{2}$ as the zero locus of the map between line bundles
\begin{displaymath}
\pi^{*} L \rightarrow \overset{2}{\wedge} S_{2}^*=\mathcal{O}_{X_{2}}(1).
\end{displaymath}
So, $Z_{2}$ is a divisor such that in $Pic(X_{2})$, $Z_{2}=u_{2}-\pi^{*}(V_{1}^{\bot})$. From the exact sequence
\begin{displaymath}
0 \rightarrow V_{1}^{\bot} \rightarrow T_{X_{1}}^* \rightarrow V_{1}^* \rightarrow 0,
\end{displaymath}
we have $c_{1}(V_{1}^{\bot})=c_{1}(T_{X_{1}}^*)-c_{1}(V_{1}^*)=c_{1}(V_{1})-c_{1}-c_{1}(T_{X_{1}/X}).$
Therefore from the exact sequence
\begin{displaymath}
0 \rightarrow T_{X_{1}/X} \rightarrow V_{1} \rightarrow S_{1} \rightarrow 0,
\end{displaymath}
we deduce that $c_{1}(V_{1}^{\bot})=c_{1}(S_{1})-c_{1}=-u_{1}-c_{1}.$ Finally we find
\begin{displaymath}
Z_{2}=u_{2}+u_{1}+c_{1}.
\end{displaymath}
\end{proof}

Now, we can prove Proposition \ref{big}.

\begin{proof}[Proof of  Proposition \ref{big}]
Let $L=F\otimes G^{-1}$ be the difference of two nef line bundles and $A$ any line bundle over a compact manifold of dimension $n$. From algebraic holomorphic Morse inequalities \cite{Tra}, for $m$ large enough, the line bundle $L^{\otimes m}\otimes A$  has a non-zero global section  as soon as $F^n-nF^{n-1}.G>0$.
Here we compute the quantity
\begin{displaymath}
F_{2}^8.Z-8F_{2}^7.G_{2}.Z=(u_{1}+5u_{2}+24h)^8.(u_{2}+u_{1}+c_{1})-8(u_{1}+5u_{2}+24h)^7.24h.(u_{2}+u_{1}+c_{1}),
\end{displaymath}
where $h$ is the class of a hyperplane section of $X\subset \bP^4$.
We apply the relations (\ref{rel1}), (\ref{rel2}), (\ref{rel3}) and Proposition \ref{chern} to express this quantity in terms of Chern classes of $X$ (the computation is made using Maple) and we get
\begin{displaymath}
F_{2}^8.Z-8F_{2}^7.G_{2}.Z=112896000\,{h}^{2}c_{{1}}+280000\,c_{{1}}c_{{2}}-70000\,c_{{3}}-
449003520\,{h}^{3}-1050000\,{c_{{1}}}^{3}.
\end{displaymath}
Since we have 
\begin{displaymath}
c_{1}=-h(d-5), c_{2}=h^2(d^2-5d+10), c_{3}=-h^3(d^3-5d^2+10d-10), h^3=d,
\end{displaymath}
(see e.g. \cite[Proposition 2]{Div}) we finally get
\begin{displaymath}
F_{2}^8.Z-8F_{2}^7.G_{2}.Z=-43246000\,{d}^{2}-2473520\,d+840000\,{d}^{4}-13300000\,{d}^{3},
\end{displaymath}
which turns out to be positive for $d\geq 19$.
\end{proof}

\subsection{Vector fields on jets spaces}
In this section, we generalize to our situation the approach of \cite{Pau} and \cite{Rou06} which represent effective versions of Siu's original ideas \cite{S04} (see also \cite{Mer} for recent generalizations). Consider $\mathcal{X}\subset \mathbb{P}^{4}\times \mathbb{P}%
^{N_{d}}$ the universal hypersurface given by the equation
\begin{equation*}
\underset{\left| \alpha \right| =d}{\sum }a_{\alpha }Z^{\alpha }=0,\text{
where }[a]\in \mathbb{P}^{N_{d}}\text{ and }[Z]\in \mathbb{P}^{4}.
\end{equation*}

We use the notations: for $\alpha =(\alpha _{0},...,\alpha _{4})\in \mathbb{N%
}^{5},$ $\left| \alpha \right| =\sum_{i}\alpha _{i}$ and if $%
Z=(Z_{0},Z_{1},...,Z_{4})$ are homogeneous coordinates on $\mathbb{P}^{4},$
then $Z^{\alpha }=\prod Z_{j}^{\alpha _{j}}.$ The variety $\mathcal{X}$ is a smooth
hypersurface of degree $(d,1)$ in $\mathbb{P}^{4}\times \mathbb{P}^{N_{d}}.$
Here, we denote by $J_{2}(\mathcal{X}):=J_{2}(\bC^2,\mathcal{X})$ the manifold of the 2-jets from $\bC^2$ to $\mathcal{X}%
, $ and $J_{2}^{v}(\mathcal{X})$ the submanifold of $J_{2}(\mathcal{X})$
consisting of 2-jets in $\mathcal{X}$ tangent to the fibers of the
projection $\pi :\mathcal{X}\rightarrow \mathbb{P}^{N_{d}}.$

Let us consider the open subset $\Omega _{0}:=(Z_{0}\neq 0)\times (a_{0d000}\neq
0)\subset \mathbb{P}^{4}\times \mathbb{P}^{N_{d}}.$ We assume that global
coordinates are given on $\mathbb{C}^{4}$ and $\mathbb{C}^{N_{d}}.$ The
equation of $\mathcal{X}$ becomes
\begin{equation*}
\mathcal{X}_{0}:=(z_{1}^{d}+\underset{\left| \alpha \right| \leq d,\text{ }%
\alpha _{1}<d}{\sum }a_{\alpha }z^{\alpha }=0).
\end{equation*}

Then the equations of $J_{2}^{v}(\mathcal{X}_{0})$ in $\mathbb{C}^{4}\times
\mathbb{C}^{N_{d}}\times \mathbb{C}^{4\times2}\times \mathbb{C}^{4\times3}$ with the coordinates $(z_{i},a_{\alpha},\xi_{j}^{(i)},\xi_{j}^{(i,l)})$ ($1\leq i\leq k \leq 2$) are given by:
\begin{equation}
\underset{\left| \alpha \right| \leq d,\text{ }a_{d000}=1}{\sum }a_{\alpha
}z^{\alpha }=0,
\end{equation}
\begin{equation}
\underset{j=1}{\overset{4}{\sum }}\underset{\left| \alpha \right| \leq d,%
\text{ }a_{d000}=1}{\sum }a_{\alpha }\frac{\partial z^{\alpha }}{\partial
z_{j}}\xi _{j}^{(i)}=0,
\end{equation}
$1\leq i \leq 2,$
\begin{equation}
\underset{j=1}{\overset{4}{\sum }}\underset{\left| \alpha \right| \leq d,%
\text{ }a_{d000}=1}{\sum }a_{\alpha }\frac{\partial z^{\alpha }}{\partial
z_{j}}\xi _{j}^{(i,l)}+\underset{j,k=1}{\overset{4}{\sum }}\underset{\left|
\alpha \right| \leq d,\text{ }a_{d000}=1}{\sum }a_{\alpha }\frac{\partial
^{2}z^{\alpha }}{\partial z_{j}\partial z_{k}}\xi _{j}^{(i)}\xi _{k}^{(l)}=0,
\end{equation}
$1\leq i\leq l \leq 2.$

Consider now a vector field
\begin{equation*}
V=\underset{\left| \alpha \right| \leq d,\text{ }\alpha _{1}<d}{\sum }%
v_{\alpha }\frac{\partial }{\partial a_{\alpha }}+\underset{j}{\sum }v_{j}%
\frac{\partial }{\partial z_{j}}+\underset{j,k}{\sum }w_{j}^{(k)}\frac{\partial
}{\partial \xi _{j}^{(k)}}+\underset{i,j,k}{\sum }w_{j}^{(i,k)}\frac{\partial
}{\partial \xi _{j}^{(i,k)}},
\end{equation*}
on the vector space  $\mathbb{C}^{4}\times
\mathbb{C}^{N_{d}}\times \mathbb{C}^{4\times2}\times \mathbb{C}^{4\times3}$.

The conditions to be
satisfied by $V$ to be tangent to $J_{2}^{v}(\mathcal{X}_{0})$ are
\begin{equation*}
\underset{\left| \alpha \right| \leq d,\text{ }\alpha _{1}<d}{\sum }%
v_{\alpha }z^{\alpha }+\underset{j=1}{\overset{4}{\sum }}\underset{\left|
\alpha \right| \leq d,\text{ }a_{d000}=1}{\sum }a_{\alpha }\frac{\partial
z^{\alpha }}{\partial z_{j}}v_{j}=0,
\end{equation*}
\begin{equation*}
\underset{j=1}{\overset{4}{\sum }}\underset{\left| \alpha \right| \leq d,%
\text{ }\alpha _{1}<d}{\sum }v_{\alpha }\frac{\partial z^{\alpha }}{\partial
z_{j}}\xi _{j}^{(i)}+\underset{j,k=1}{\overset{4}{\sum }}\underset{\left|
\alpha \right| \leq d,\text{ }a_{d000}=1}{\sum }a_{\alpha }\frac{\partial
^{2}z^{\alpha }}{\partial z_{j}\partial z_{k}}v_{j}\xi _{k}^{(i)}+\underset{%
j=1}{\overset{4}{\sum }}\underset{\left| \alpha \right| \leq d,\text{ }%
a_{d000}=1}{\sum }a_{\alpha }\frac{\partial z^{\alpha }}{\partial z_{j}}%
w_{j}^{(i)}=0,
\end{equation*}
$1\leq i \leq 2$, and
\begin{eqnarray*}
\underset{\left| \alpha \right| \leq d,\text{ }\alpha _{1}<d}{\sum }(%
\underset{j=1}{\overset{4}{\sum }}\frac{\partial z^{\alpha }}{\partial z_{j}}%
\xi _{j}^{(i,l)}+\underset{j,k=1}{\overset{4}{\sum }}\frac{\partial
^{2}z^{\alpha }}{\partial z_{j}\partial z_{k}}\xi _{j}^{(i)}\xi
_{k}^{(l)})v_{\alpha } && \\
+\underset{j=1}{\overset{4}{\sum }}\underset{\left| \alpha \right| \leq d,%
\text{ }a_{d000}=1}{\sum }a_{\alpha }(\underset{k=1}{\overset{4}{\sum }}%
\frac{\partial ^{2}z^{\alpha }}{\partial z_{j}\partial z_{k}}\xi _{k}^{(i,l)}+%
\underset{k,m=1}{\overset{4}{\sum }}\frac{\partial ^{3}z^{\alpha }}{\partial
z_{j}\partial z_{k}\partial z_{l}}\xi _{k}^{(i)}\xi _{m}^{(l)})v_{j} && \\
+\underset{\left| \alpha \right| \leq d,\text{ }a_{d000}=1}{\sum }(\underset{%
j,k=1}{\overset{4}{\sum }}a_{\alpha }\frac{\partial ^{2}z^{\alpha }}{%
\partial z_{j}\partial z_{k}}(w_{j}^{(i)}\xi _{k}^{(l)}+w_{k}^{(i)}\xi
_{j}^{(l)})+\underset{j=1}{\overset{4}{\sum }}a_{\alpha }\frac{\partial
z^{\alpha }}{\partial z_{j}}w_{j}^{(i,l)}) &=&0,
\end{eqnarray*}
$1 \leq i \leq l \leq 2.$

\bigskip
We obtain exactly the same type of equations as in \cite{Pau}. Therefore using the computations of \cite{Pau}, we obtain meromorphic vector fields on $T_{J_{2}^{v}(\mathcal{X})}.$ We recall the main facts for the convenience of the reader and refer to the above mentioned article for the details. 

The first package of vector fields tangent to $J_{2}^{v}(\mathcal{X})$ is the following.

For $\alpha _{1}\geq 3:$%
\begin{equation*}
V_{\alpha }^{300}:=\frac{\partial }{\partial a_{\alpha }}-3z_{1}\frac{%
\partial }{\partial a_{\alpha -\delta _{1}}}+3z_{1}^{2}\frac{\partial }{%
\partial a_{\alpha -2\delta _{1}}}-z_{1}^{3}\frac{\partial }{\partial
a_{\alpha -3\delta _{1}}}.
\end{equation*}

For $\alpha _{1}\geq 2,\alpha _{2}\geq 1:$%
\begin{eqnarray*}
V_{\alpha }^{210} &:&=\frac{\partial }{\partial a_{\alpha }}-2z_{1}\frac{%
\partial }{\partial a_{\alpha -\delta _{1}}}-z_{2}\frac{\partial }{\partial
a_{\alpha -\delta _{2}}}+ \\
&&+z_{1}^{2}\frac{\partial }{\partial a_{\alpha -2\delta _{1}}}+2z_{1}z_{2}%
\frac{\partial }{\partial a_{\alpha -\delta _{1}-\delta _{2}}}-z_{1}^{2}z_{2}%
\frac{\partial }{\partial a_{\alpha -2\delta _{1}-\delta _{2}}}.
\end{eqnarray*}

For $\alpha _{1}\geq 1,\alpha _{2}\geq 1,\alpha _{3}\geq 1:$%
\begin{eqnarray*}
V_{\alpha }^{111} &:&=\frac{\partial }{\partial a_{\alpha }}-z_{1}\frac{%
\partial }{\partial a_{\alpha -\delta _{1}}}-z_{2}\frac{\partial }{\partial
a_{\alpha -\delta _{2}}}-z_{3}\frac{\partial }{\partial a_{\alpha -\delta
_{3}}} \\
&&+z_{1}z_{2}\frac{\partial }{\partial a_{\alpha -\delta _{1}-\delta _{2}}}%
+z_{1}z_{3}\frac{\partial }{\partial a_{\alpha -\delta _{1}-\delta _{3}}} \\
&&+z_{2}z_{3}\frac{\partial }{\partial a_{\alpha -\delta _{2}-\delta _{3}}}%
-z_{1}z_{2}z_{3}\frac{\partial }{\partial a_{\alpha -\delta _{1}-\delta
_{2}-\delta _{3}}}.
\end{eqnarray*}

We obtain similar vector fields by permuting the $z_{i}$ and changing $\alpha $ as indicated by the permutation.
The order of the pole is equal to 3.

Another family is obtained in the following way.
Consider a $4\times 4$-matrix $A=(A_{j}^{k})\in {M}_{4}(%
\mathbb{C})$ and let $\widetilde{V}:=\underset{j}{\sum }w_{j}^{(k)}\frac{\partial
}{\partial \xi _{j}^{(k)}}+\underset{j}{\sum }w_{j}^{(i,k)}\frac{\partial
}{\partial \xi _{j}^{(i,k)}},$ where $w^{(k)}:=A\xi ^{(k)}$ and $w^{(i,k)}:=A\xi ^{(i,k)}.$

\begin{lemma}[lemma 1.1 of \cite{Pau}]
There exist polynomials $v_{\alpha }(z,a):=\underset{\left| \beta \right|
\leq 3}{\sum }v_{\beta }^{\alpha }(a)z^{\beta }$ where each coefficient $%
v_{\beta }^{\alpha }$ has degree at most 1 in the variables $(a_{\gamma })$
such that
\begin{equation*}
V:=\underset{\alpha }{\sum }v_{\alpha }(z,a)\frac{\partial }{\partial
a_{\alpha }}+\widetilde{V}
\end{equation*}
is tangent to $J_{2}^{v}(\mathcal{X}_{0})$ at each point.
\end{lemma}

We define the affine algebraic set
\begin{displaymath}
\Sigma_{0}:=\underset{1\leq i,k \leq 2} \bigcap \Sigma _{i,k} \subset J_{2}^{v}(\mathcal{X}),
\end{displaymath}
where
\begin{displaymath}
\Sigma _{i,k}:=\{(z,a,\xi ^{(i)},\xi ^{(i,k)})\in
J_{2}^{v}(\mathcal{X}) /\xi ^{(i)}\wedge \xi^{(i,k)}=0\}, 1\leq i,k \leq 2.
\end{displaymath}

We remark that we have
\begin{lemma}\label{rk}
If the lifting of the map of maximal rank $f:\bC^2 \rightarrow X_{a}\subset \mathcal{X}$ lies inside $\Sigma_{0}$, then the image of $f$ is contained in a proper suvariety of ${X}_a$.
\end{lemma}

\begin{proof}
By hypothesis $f$ verifies the equations
\begin{equation}
\label{col}
\frac{\partial{f}}{\partial{t_{i}}}\wedge \frac{\partial^2{f}}{\partial{t_{i}}\partial{t_{k}}}=0,
\end{equation}
$1\leq i,k \leq 2$.
Since $f$ is of maximal rank, this implies that $\frac{\partial^2{f}}{\partial{t_{1}}\partial{t_{2}}}=0$ and $f=g(t_1)+h(t_{2}).$ From (\ref{col}) we get
\begin{eqnarray*}
\frac{\partial{g}}{\partial{t_{1}}}\wedge \frac{\partial^2{g}}{\partial{t_{1}}^2}=0,\\
\frac{\partial{h}}{\partial{t_{2}}}\wedge \frac{\partial^2{h}}{\partial{t_{2}}^2}=0.
\end{eqnarray*}
This provides the algebraic degeneracy of $f$.
\end{proof}

Now, we have the following.
\begin{proposition}\label{gg}
 The vector bundle $T_{J_{2}^{v}(\mathcal{X)}}\otimes
\mathcal{O}_{\mathbb{P}^{4}}(7)\otimes \mathcal{O}_{\mathbb{P}%
^{N_{d}}}(q )$ is generated by its global sections on $J_{2}^{v}(\mathcal{%
X})\backslash \Sigma ,$ where $\Sigma $ is the closure of $\Sigma _{0}$ and $q$ a sufficiently large integer.
\end{proposition}

\begin{proof}
From the previous lemmas, we are reduced to consider vector fields of the form $V=\underset{\left|
\alpha \right| \leq 2}{\sum }v_{\alpha }\frac{\partial }{\partial a_{\alpha }%
}.$ The conditions for $V$ to be tangent to $J_{2}^{v}(\mathcal{X}_{0})$ are
\begin{equation*}
\underset{\left| \alpha \right| \leq 2}{\sum }v_{\alpha }z^{\alpha }=0,
\end{equation*}
\begin{equation*}
\underset{j=1}{\overset{4}{\sum }}\underset{\left| \alpha \right| \leq d,%
\text{ }\alpha _{1}<d}{\sum }v_{\alpha }\frac{\partial z^{\alpha }}{\partial
z_{j}}\xi _{j}^{(i)}=0,
\end{equation*}
and
\begin{equation*}
\underset{\left| \alpha \right| \leq 2}{\sum }(\underset{j=1}{\overset{4}{%
\sum }}\frac{\partial z^{\alpha }}{\partial z_{j}}\xi _{j}^{(i,l)}+\underset{%
j,k=1}{\overset{4}{\sum }}\frac{\partial ^{2}z^{\alpha }}{\partial
z_{j}\partial z_{k}}\xi _{j}^{(i)}\xi _{k}^{(l)})v_{\alpha }=0.
\end{equation*}

We are interested in the global generation outside $\Sigma$ so we can assume that for some $1\leq i,l \leq 2$, we have $W_{12}:=\xi _{1}^{(i)}\xi _{2}^{(i,l)}-\xi _{1}^{(i,l)}\xi_{2}^{(i)}\neq 0$. Then we can solve the previous system with $v_{000}, v_{100}$ and $v_{010}$ as unknowns. By the Cramer rule, each
of the previous quantity is a linear combination of the 
$v_{\alpha },$ $\left| \alpha \right| \leq 2,$ $\alpha \neq (000), (100), (010)$ 
with coefficients rational functions in $z,\xi ^{(i)},\xi ^{(i,l)}.$ The denominator is $W_{12}$ and the numerator is a polynomial whose monomials have either degree at most $2$ in $z$ and at most one in $\xi ^{(i)}$ and $\xi ^{(i,l)}$, or degree $1$ in $z$ and three in $\xi ^{(i)}$ variables. Thus we get order $7$ for the above vector field.
\end{proof}

\subsection{Proof of Theorem \ref{thm:deg2}}
First we need a lemma which is a slight modification of Theorem \ref{jdiff}.

\begin{lemma}
Let $X\subset \bP^4$ be an hypersuface of degree $d\geq 6$. Then there exists non-zero global sections in
$$H^0(Z_{2},\mathcal{O}_{Z_{2}}(6m)\otimes\pi_{2,0}^*K_{X}^{-\delta m}) \subset H^0(X,E_{2,2,6m}\otimes K_{X}^{-\delta m}),$$ for $m$ large enough if
\begin{eqnarray}\nonumber
\alpha(d,\delta)=-742000\,{d}^{4}\delta-36400\,d\delta+1862000\,{d}^{3}\delta-43246000\,{d}^{2}\\\nonumber
-13300000\,{d}^{3}+840000\,{d}^{4}+9247280\,{d}^{2}\delta-2473520\,d >0.
\end{eqnarray}
\end{lemma}
\begin{proof}
We apply algebraic holomorphic Morse inequalities with $F_{2}=\cO_{X^2_2}(5,1)\otimes \pi_{0,2}^* \cO_X(24)$, $G_{2}=\pi_{2,0}^*(\cO_{X}(24)+K_{X}^{\delta})$ and $L_{2}=F_{2}\otimes G_{2}^{-1}$ on the effective locus $Z_{2}$. To finish the proof, we compute everything as in Proposition \ref{big} and use the injection (\ref{injection}).
\end{proof}

Now, we can finish the proof of Theorem  \ref{thm:deg2}.

\begin{proof}[Proof of Theorem  \ref{thm:deg2}]
Let us consider a maximal rank holomorphic map $f:\mathbb{C}^2\rightarrow X$ in a very generic
hypersurface of $\mathbb{P}^{4}$ of degree $d \geq 93.$

\bigskip
We have a section
\begin{equation*}
P \in H^{0}(Z_{2},\mathcal{O}_{Z_{2}}(6m)\otimes \pi _{2}^{\ast }K_{X}^{-\delta m}) \subset H^{0}(X,E_{2,2,6m}\otimes K_{X}^{-\delta m}),
\end{equation*}
with zero set $Y$ and vanishing order $\delta m(d-5).$ Consider the family
\begin{equation*}
\mathcal{X}\subset \mathbb{P}^{4}\times \mathbb{P}^{N_{d}}
\end{equation*}
of hypersurfaces of degree $d$ in $\mathbb{P}^{4}.$ 
We denote by $X_a$ the fiber of $\mathcal{X}$  over the parameter $a\in \mathbb{P}^{N_{d}}
$.
General semicontinuity
arguments concerning the cohomology groups show the existence of a Zariski
open set $U_{d}\subset \mathbb{P}^{N_{d}}$ such that for any $a\in U_{d},$
there exists a divisor
\begin{equation*}
Y_{a}=(P_{a}=0)\subset ({Z}_{a})_{2},
\end{equation*}
where
$(Z_a)$ is the effective locus inside $(X_a)^2_2$ and $(P_{a})_{a\in U_{d}}$ is a family of sections
\begin{equation*}
P_{a}\in H^{0}(({Z}_{a})_{2},\mathcal{O}_{({Z}%
_{a})_{2}}(6m)\otimes \pi _{2}^{\ast }K_{({X}_{a})}^{-\delta m}),
\end{equation*} 
varying holomorphically over $U_d$. We
consider $P$ as a holomorphic function on $J_{2}({X}_{a}).$ The
vanishing order of this function as a function of $%
\xi^{(i)},\xi^{(i,l)}$ $(1\leq i,l\leq 2)$ is no more than $12m$ at a
generic point of $\mathcal{X}_{a}.$ We have $\Im(f_{[2]}) \subset
Z_{a}.$ If $\Im(f_{[2]})$ lies in $\Sigma$, $f$ is algebraically degenerated
as we saw in Lemma \ref{rk}. So, we can suppose it is not the case.

If we take a vector field $v$ given by Proposition \ref{gg}, we can differentiate $P$. We see that $vP$ corresponds to an invariant
differential operator and its restriction to $(P_{a}=0)$ can be
seen as a section of the bundle
\begin{equation*}
\mathcal{O}_{({Z}_{a})_{2}}(6m)\otimes \mathcal{O}_{\mathbb{P}%
^{4}}(7-\delta m(d-5))_{|Y_{a}}.
\end{equation*}
If $\Im(f_{[2]})$ does not lie in the singular locus of $(P_{a}=0)$, we invoke Proposition \ref{gg} which gives the global generation of
\begin{equation*}
T_{J_{2}^{v}(\mathcal{X)}}\otimes \mathcal{O}_{\mathbb{P}^{4}}(7)\otimes
\mathcal{O}_{\mathbb{P}^{N_{d}}}(q)
\end{equation*}
on $J_{2}^{v}(\mathcal{X})\backslash \Sigma$, to obtain a vector field $v$ such that $vP$ is non zero at some point of $\Im(f_{[2]})$. If
$$7-\delta m(d-5)<0,$$ 
we obtain a contradiction from Theorem \ref{thm:deg} because $vP$ should be zero at this point.

If $\Im(f_{[2]})$ lies in the singular locus of $(P_{a}=0)$, we have to differentiate more, but not more than $12m$. We can find global meromorphic vector fields $%
v_{1},...,v_{q}$ $(q\leq 12m)$ and differentiate $P$ with respect to these vector fields
such that $v_{1}...v_{q}P$ is not zero at some point of $\Im(f_{[2]})$. Assume that the vanishing order of $P$ is larger than the sum of the pole
order of the $v_{i}$ in the fiber direction of $\pi :\mathcal{X}\rightarrow
\mathbb{P}^{N_{d}}.$ Then as before the restriction of $v_{1}...v_{q}P$  can be seen as the restriction of a global section of  $\mathcal{O}_{(\mathcal{Z}_{a})_{2}}(6m)$ which vanishes on the pull-back of an ample divisor.
Therefore the lifting $f_{[2]}$ should be in its zero set. Thus $f_{[2]}$ should be in the zero section of $J_{2}(\mathcal{X}_{a})$ over
a generic point of $\mathcal{X}_{a}.$

To obtain the algebraic degeneracy of $f$, we just have to see when the vanishing order of $P$ is
larger than the sum of the pole order of the $v_{i}.$ This will be verified
if
$$
7\times 12m-\delta m(d-5)<0,
$$
i.e.
\begin{equation*}
\delta (d-5)>84.
\end{equation*}
So we want $\delta >\frac{84}{(d-5)}$ and $\alpha (d,\delta )>0.$ This is
the case for $d\geq 93.$

To finish the proof we remark now that the Zariski closure of $f$ is a surface $S$ which cannot be of general type by a classical result \cite{KO75} of Kobayashi and Ochiai. But by a result of Ein \cite{Ein} (later sharpened in \cite{V'}, \cite{P}, and \cite{CR})  all subvarieties of a very general hypersurface of degree $\geq 2n$ in $\bP^n$ are of general type, for all $n\geq 3$. So we obtain a contradiction and there is no maximal rank map $f:\mathbb{C}^2\rightarrow X$ in a very generic
hypersurface of $\mathbb{P}^{4}$ of degree $d \geq 93.$
\end{proof}

\bigskip

\noindent
{\tt pacienza@math.u-strasbg.fr} and 
{\tt rousseau@math.u-strasbg.fr}\\
Institut de Recherche Math\'ematique Avanc\'ee\\
Universit\'e L. Pasteur et CNRS\\
7, rue Ren\'e Descartes, 67084 Strasbourg C\'edex - FRANCE

\end{document}